\documentclass[12pt,english]{elsarticle}
\usepackage[utf8]{inputenc}
\usepackage{amsmath,amssymb,multirow}
\usepackage{color}
\usepackage{psfrag}
\usepackage{soul}

\usepackage{amssymb}

\renewcommand\appendix{\par
  \setcounter{section}{0}
  \setcounter{subsection}{0}
  \setcounter{figure}{0}
  \setcounter{table}{0}
  \renewcommand\thesection{Appendix \Alph{section}}
  \renewcommand\thefigure{\Alph{section}\arabic{figure}}
  \renewcommand\thetable{\Alph{section}\arabic{table}}
}

\begin{document}
\graphicspath{{FIG/}}

\begin{frontmatter}



\title{An adaptive multiresolution method for ideal magnetohydrodynamics using divergence cleaning
with parabolic-hyperbolic correction} 


\author[CAP,CTE,INPE]{Anna Karina Fontes Gomes}
\ead{annakfg@gmail.com}
\author[LAC,CTE,INPE]{Margarete Oliveira Domingues}
\ead{margarete@lac.inpe.br, margarete.oliveira.domingues@gmail.com}
%
\author[MP2]{\\ Kai Schneider}
\ead{kschneid@cmi.univ-mrs.fr}
\author[DGE,CEA,INPE]{Odim Mendes}
\ead{odim@dge.inpe.br,o-mendes@hotmail.com}
\author[DLR]{Ralf Deiterding}
\ead{ralf.deiterding@dlr.de}
\address[CAP]{Pós-Graduação em Computação Aplicada (CAP)}
\address[LAC]{Laboratório Associado de Computação e Matemática Aplicada (LAC)}
\address[CTE]{Coordenadoria dos Laboratórios Associados (CTE)}
\address[DGE]{Divisão de Geofísica Espacial,
Coordenação de Ciências Espaciais(CEA)}
\address[INPE]{Instituto Nacional de Pesquisas Espaciais (INPE), 
Av. dos Astronautas 1758, 12227-010 São José dos Campos, São Paulo, Brazil}
\address[MP2]{M2P2--CNRS \& Centre de Math\'ematiques et d'Informatique (CMI), Aix-Marseille Universit\'e, 38 rue F. Joliot--Curie, 13451 Marseille Cedex 20, France}
\address[DLR]{German Aerospace Center (DLR), Institute of Aerodynamics and Flow Technology, Bunsenstr. 10, 37073 Göttingen, Germany}

%
%

\begin{abstract}
We present an adaptive multiresolution method for the numerical simulation of ideal
magnetohydrodynamics in two space dimensions. The discretization uses a finite volume
scheme based on a Cartesian mesh and an explicit compact Runge–Kutta scheme for
time integration. Harten’s cell average multiresolution allows to introduce a locally
refined spatial mesh while controlling the error. The incompressibility of the magnetic
field is controlled by using a Generalized Lagrangian Multiplier (GLM) approach with a
mixed hyperbolic–parabolic correction. Different applications to two-dimensional problems
illustrate the properties of the method. For each application CPU time and memory savings
are reported and numerical aspects of the method are discussed. The accuracy of the
adaptive computations is assessed by comparison with reference solutions computed on a
regular fine mesh.
\end{abstract}

\end{frontmatter}



\section{Introduction}
The magnetohydrodynamic (MHD) equations, which consist of the compressible Euler equations of hydrodynamics coupled with the Maxwell equations of electrodynamics, are used for
 mathematical modeling of numerous phenomena encountered in our daily life.
Prominent examples can be found in the physics of the Sun-Earth's electrodynamical interaction chain,
and in the dynamo action caused by motion of liquid metal inside the mantle of the Earth, which generates its magnetic field.
The numerical challenge for solving the ideal MHD equations, a coupled set of nonlinear Partial Differential Equations (PDEs),
is the presence of multiple spatial and temporal scales.
The complex character of boundary conditions of the magnetic field, in comparison to that one for the classical hydrodynamics, requires even more sophisticated approaches.
In a surrounding vacuum, for example, the magnetic field does not vanish, it only decays.
Thus, at the boundary it has to be matched with the field of the fluid region.
A second difficulty is to maintain  the incompressibility of the magnetic field numerically, which is imposed by Gauss' law.
Therefore, in the numerical simulations, special attention has to be paid to this incompressibility, because, as shown in practice, uncontrolled divergence errors can modify the underlying physics. 
For details we refer the reader to, \textit{e.g.}, \cite{Balsara:2009,BrBa80, Powell:1999,Toth:2000JCP}.
Typically, projection methods based on the Helmholtz decomposition are used.
These methods are computationally demanding, especially in three-dimension,  because  the solution of  an elliptic problem requires a Poisson equation solver.
An alternative method is the divergence cleaning one, which is based on Lagrangian multipliers.
In the finite element context, Assous et al. \cite{Assous1993222} introduced this approach for time-dependent Maxwell equations.
Several variants can be found in the literature \cite{Balsara:2009,BrBa80,Powell:1994,Toth:2000JCP}.

In the current paper we  apply the multiresolution approach to an ideal MHD numerical model called the Generalized
Lagrange Multiplier (GLM) with a mixed hyperbolic-parabolic correction proposed by Dedner et al. \cite{Dedneretal:2002}
%
to deal with the magnetic field incompressibility condition.
The ideas of the Lagrangian multiplier formulation in this context were introduced by Munz et al. \cite{Munzetal:2000} in the context of Maxwell equations. 
With the motivation to reduce CPU time and memory requirements, we
use an auto-adaptive discretization which is based on the multiresolution representation. 
The underlying time dependent conservation laws are discretized with finite volume schemes and local grid refinement is triggered by multiresolution analysis of the cell averages and thresholding of the resulting coefficients. 
The adaptive refinement/mesh tracks steep gradients in the solution of the equation and allows automatic error control.
For reviews on multiresolution techniques for PDEs we refer to \cite{Harten:1995,Harten:1996,Mueller:2003,DGRSESAIM:2011} and references therein.

Preliminary results for a quasi-one dimensional MHD Riemann problem with exact solution have been presented in \cite{Dominguesetal:2013}, which showed the feasibility of using adaptive discretizations and magnetic field divergence cleaning for extended GLM--MHD with local and controlled time methods. In its extended form, source terms similar to those in \cite{Powell:1994} are introduced.
The starting point is the adaptive multiresolution code originally developed by Roussel et al. \cite{RSTB03} in which the Maxwell equations governing the magnetic field have been included \cite{Gomes:2012:AnMuAd}.
In the present work, we have chosen the GLM--MHD approach instead of its extended version, because the divergence errors and the solution obtained for both cases are almost the same for the studied problem.
A similar choice is suggested in the conclusion in \cite{Dedneretal:2002}.
The resulting new method has been applied to a two-dimensional Riemann test problem, for which a reference solution on a fine grid has been computed. The accuracy of the adaptive computations has been assessed and their efficiency in terms of memory compression compared to a finite volume scheme on a regular grid has been analyzed.

The paper is organized as follows: After a presentation of the governing ideal MHD equations in Section~\ref{sec:MHD}, we recall the divergence cleaning technique based on the GLM formulation in Section~\ref{sec:GLM}. In Section~\ref{sec:SpaceTime} space and time discretizations are briefly described
together with the GLM discretization. In Section~\ref{sec:numerical}, numerical
results are presented. In the last section, some conclusions are drawn and perspectives for future work are presented.

\section{Governing equations}
\label{sec:MHD}
The ideal magnetohydrodynamics equations describe the dynamics of a compressible, inviscid and perfectly electrically conducting fluid interacting with a magnetic field, see, \textit{e.g.} \cite{freidberg2014ideal}.   
The equations combine the Euler equations with the Maxwell equations. 
The latter yields an evolution equation for the magnetic field, also called induction equation, and an incompressibility constraint using Gauss' law. 
The system of MHD equations is given by
\begin{small}
\begin{subequations}
\label{eq:MHD}
\begin{eqnarray}
\displaystyle\frac{\partial\rho}{\partial t} +\nabla\cdot (\rho {\bf u}) = 0,& 
 \text{(Mass conservation)} \label{eq:mass}
\\
\displaystyle\frac{\partial E}{\partial t} +\nabla\cdot\left[\left(E + p + \frac{{\bf B \cdot B}}{2}\right){\bf u} - \left({\bf u \cdot B} \right) {\bf B}\right]= 0,& \text{(Energy conservation)}\label{eq:energy}
\\
\displaystyle\frac{\partial\rho {\bf u}}{\partial t} 
+ \nabla\cdot\left[\rho{\bf u^t u} + \left(p + \frac{{\bf B \cdot B}}{2}\right) {\bf I}  - {\bf B^t B}\right]
= {\bf 0},& 
 \text{{\footnotesize (Momentum conservation)}} \label{eq:momentum}
\\
\displaystyle\frac{\partial{\bf B}}{\partial t} 
+ \nabla\cdot\left({\bf u^t B-B^t u}\right)={\bf 0},&
 \text{(Induction equation)} \label{eq:induction}
		\end{eqnarray}
	\end{subequations}
\end{small}
\noindent 
where $\rho$ represents   density, $p$ the pressure, ${\bf u} =(u_x, u_y, u_z)$  the velocity vector, ${\bf B} = (B_x,B_y,B_z)$ the magnetic field vector, and $t$ denotes the transposition. The identity tensor of order 2 is denoted by ${\bf I}$ (the unit dyadid, that here corresponds to the unit matrix $3\times 3$), and  $\gamma$  the adiabatic constant ($\gamma>1$). 
The pressure is given by the constitutive law $
p = \left(\gamma-1 \right) 
\left( 
E -\rho\frac{{\bf u\cdot u}}{2}-\frac{{\bf B\cdot B}}{2} 
\right).
\label{eq:pressure}
$
The above system is completed by suitable initial and boundary conditions.
In this paper this system is considered in its two-dimensional form, \textit{i.e.}, the quantities depend on two variables only ($x$ and $y$).

In this classical MHD model, the magnetic field has to satisfy the divergence constraint 
\begin{equation}
\nabla\cdot{\bf B}={\bf 0}. 
\label{eq:divconstraint}
\end{equation}	  
which implies the non-existence of magnetic monopoles. By rewriting the induction equation, we have $\displaystyle\frac{\partial{\bf B}}{\partial t} + \nabla \times\left({\bf B\times u}\right)={\bf 0}$. 
Therefore, the  application of  the divergence operator yields  $\displaystyle\frac{\partial}{\partial t}  \left(\nabla \cdot \bf B \right) = 0$, as $ \nabla \cdot \left( \nabla \times  \;  \right) \equiv 0$. 
This formulation shows that if the initial condition of the magnetic field is divergence-free, the system  will remain divergence-free along the evolution.
However, numerically the incompressibility of the magnetic field is not necessarily preserved, and thus,  non-physical results could be obtained or the computations may even become unstable \cite{BrBa80}.
Since the 1980ies typical numerical MHD methodologies consider the enforcement of the divergence-free constraint. 
There are many techniques to perform the divergence cleaning in the MHD numerical models \cite{Toth2012870}.
In the context of this study, we have in mind the application of the multiresolution method based on a finite volume discretization with explicit time integration. Thus, the  technique developed in Dedner et al.\cite{Dedneretal:2002} called GLM--MHD with the mixed parabolic-hyperbolic correction, is well suited. Details are given in the next section.%

\section{Generalized Lagrangian multipliers for divergence cleaning}
\label{sec:GLM}

Dedner et al. \cite{Dedneretal:2002} proposed the GLM formulation with the hyperbolic-parabolic correction.
Its implementation into a pre-existing MHD model is straightforward. An additional scalar field $\psi$ is introduced, which couples the divergence constraint equation (Eq.~\ref{eq:divconstraint}) to Faraday's law, modifying the induction equation (Eq.~\ref{eq:induction}).
Moreover, some source terms are added similarly to what was proposed in \cite{Powell:1994}.
The model contains one parameter related to the hyperbolic correction, namely $c_h$, responsible for the propagation of the divergence errors, and another one related to the parabolic correction $c_p$, responsible for the damping of the monopoles.
The remaining terms in the equations remain unchanged. 
The conservative characteristic of this system is not lost for the GLM approach.

The resulting GLM--MHD equations written in two-dimensional form read
\begin{subequations}
\label{glm}
\begin{small}
\begin{eqnarray}
&& \displaystyle\frac{\partial\rho}{\partial t}
+\frac{\partial\rho u_x}{\partial x}+\frac{\partial \rho u_y}{\partial y}=0,
\label{eq:eglmMass}
\\
&& \displaystyle\frac{\partial E}{\partial t} +  \frac{\partial}{\partial x}\left[\left(E + p + \frac{{\bf B \cdot B}}{2}\right){\bf u_x} - \left({\bf u \cdot B} \right) { B_x}\right] + \nonumber\\ &&\qquad \;\;\;\frac{\partial\rho}{\partial y}\left[\left(E + p + \frac{{\bf B \cdot B}}{2}\right){\bf u_y} - \left({\bf u \cdot B} \right) { B_y}\right] = 0, 
\label{eq:eglmEnergy}
\\
&& \displaystyle\frac{\partial \left(\rho { u_x}\right)}{\partial t} + \frac{\partial}{\partial x}\left[\rho u_x^2 + p \left( p + \frac{{\bf B\cdot B}}{2} \right) \!-\! B_x^2\right] + \frac{\partial}{\partial y}\left(\rho u_x u_y - B_xB_y\right)  = { 0}, 
\label{eq:eglmMomentumX}
\\
&& \displaystyle\frac{\partial \left(\rho { u_y}\right)}{\partial t} + \frac{\partial}{\partial x}\left(\rho u_x u_y \!-\! B_xB_y\right)  + \frac{\partial}{\partial y}\left[\rho u_y^2 + p \left( p + \frac{{\bf B\cdot B}}{2} \right)- B_y^2\right] = { 0}, 
\label{eq:eglmMomentumY}
\\
&& \displaystyle\frac{\partial \left(\rho { u_z}\right)}{\partial t} + \frac{\partial}{\partial x}\left(\rho u_z u_x - B_zB_x\right) + \frac{\partial}{\partial y}\left(\rho u_z u_y - B_zB_y\right)  = { 0}, 
\label{eq:eglmMomentumZ}
\\
&& \displaystyle\frac{\partial{ B_x}}{\partial t} + \frac{\partial \psi}{\partial x} + \frac{\partial}{\partial y}\left( u_yB_x - B_y u_x \right)={ 0}, 
\label{eq:eglmInductionX}
\\
&& \displaystyle\frac{\partial{ B_y}}{\partial t} + \frac{\partial}{\partial x}\left( u_xB_y - B_x u_y \right)+\frac{\partial \psi}{\partial y} ={ 0}, 
\label{eq:eglmInductionY}
\\
&& \displaystyle\frac{\partial{ B_z}}{\partial t} + \frac{\partial}{\partial x}\left( u_xB_z - B_z u_x \right) + \frac{\partial}{\partial y}\left( u_yB_z - B_y u_z \right)   ={ 0}, 
\label{eq:eglmInductionZ}
\\
&& \frac{\partial\psi}{\partial t}+c_h^2\left(\frac{\partial B_x}{\partial x} + \frac{\partial B_y}{\partial y}\right)=-\frac{c_h^2}{c_p^2}\psi,
\label{eq:eglmDivConstraintHP}
\end{eqnarray}
\end{small}
\end{subequations}
\noindent where ${\bf B \cdot \bf B}=B_x^2 + B_y^2 + B_z^2$, ${\bf u \cdot \bf B} = u_xB_x + u_yB_y + u_zB_z$, $c_p$ and $c_h$ are the parabolic-hyperbolic parameters, with $c_h >0$.
In \cite{Dedneretal:2002} it is defined as 
\[
c_h = c_h(t) := c_{CFL} \frac{\min\{\Delta x,\Delta y\}}{\Delta t},
\label{eq:ch}
\]
where $c_{CFL} \in (0,1)$, $\Delta x$ and $\Delta y$ are the space step in $x-$ and $y-$direction, respectively, $\Delta t$ is the time step. If the parameter $c_h$ is defined, as for instance in Eq.~\ref{eq:ch}, then $c_p$ is a free parameter in Eq.~\ref{eq:eglmDivConstraintHP}. We follow a choice proposed in \cite{Dedneretal:2002} to avoid that $c_p$ is strongly dependent on the mesh size and the scheme used. Their numerical experiments showed that choosing $c_p^2/c_h=0.18$, mirrors properly the ratio between hyperbolic and parabolic effects. With this choice in the one-dimensional case the damping of the divergence errors occurs on the time scale $c_p\sqrt{t}$ and the transport of the divergence errors to the boundary takes place on the time scale $c_h t$ (as discussed in \cite{Dedneretal:2002}, Appendix A.16 and A.19). However, other possible choices of these parameters can be found in \cite{Tricco20127214,Dedneretal:2003} and for the CTU--GLM approach in \cite{mignone2010second}.

Considering the vector of conservative quantities
${\bf Q}= (\rho, E, \rho {\bf u},{\bf B}, \psi)$, the GLM--MHD system could be written compactly as
\[
\displaystyle \frac{\partial \bf{Q}}{\partial t} 
+ \nabla \cdot {\bf{F}}({\bf Q}) = {\bf S}({\bf Q}),  
\]
where $\bf F(Q)$ is the physical flux, and {\bf S}({\bf Q}) contains all source terms.

\section{Adaptive space and time discretization}
\label{sec:SpaceTime}
A finite volume discretization of the GLM--MHD system is applied, which  results in a system of ordinary differential equations. 
Approximate solutions at a sequence of time instants $t^{n}$ are obtained by using an explicit ordinary differential equation solver. 
Here, an explicit Runge-Kutta scheme of second order is used.

In the GLM--MHD Finite Volume (FV) reference scheme, we consider the initial value of the variable $\psi$ as zero. The parameter $c_h$ has a strong influence in the correction.
In each time step, we compute the parameter $c_h$, then the GLM--MHD system is solved. First, a dimensional splitting is performed in $x$-direction, where the fluxes in the interface are treated and the solution updated. This procedure follows the steps:
\begin{enumerate}
\item The component of the magnetic field $B_x$ in the $x$-direction flux (Eq.~\ref{eq:eglmInductionX}), and the divergence constraint equation (Eq.~\ref{eq:eglmDivConstraintHP}), are decoupled from the other variables.
These two equations form the system 
\begin{eqnarray}
\displaystyle\frac{\partial{ B_x}}{\partial t} + \frac{\partial \psi}{\partial x} &=&{ 0}, 
\label{eq:eglmInductionX1D}
\\
\frac{\partial\psi}{\partial t}+c_h^2\frac{\partial B_x}{\partial x}&=&-\frac{c_h^2}{c_p^2}\psi,
\end{eqnarray}
such that the local Riemann problem can be solved analytically, where the numerical flux in the interface is $(\psi_m, c_h^2 B_{x,m})$ for $B_x$ and $\psi$. Similarly as what is described in \cite{Dedneretal:2002}, we have
\begin{equation}
\left( \begin{array}{c} 
        B_{x,m}\\ 
        \psi_m
        \end{array}
\right)=
\left( \begin{array}{c} 
        B_{x,L}\\ 
        \psi_L
        \end{array}
\right) + 
\left( \begin{array}{c} 
        \frac{1}{2}(B_{x,R} - B_{x,L}) - \frac{1}{2c_h}(\psi_R - \psi_L)\\ 
         \frac{1}{2}(\psi_R - \psi_L) - \frac{c_h}{2}(B_{x,R} - B_{x,L})
        \end{array}
\right)        
\end{equation}
where the sub-index $L,R$ are related to the left or right-hand state.
\\ \\
\item  Therefore, the numerical flux is evaluated in two steps. First we compute the numerical flux not considering the $B_x$ and $\psi$ equations as described above, then we add the numerical flux in the interface. 
In this work, we use the  Harten-Lax-van Leer-Discontinuities numerical flux (HLLD) with four intermediary states ${\bf Q}_L^{\star}$, ${\bf Q}_L^{\star\star}$, ${\bf Q}_R^{\star\star}$ and ${\bf Q}_R^{\star}$, divided by the waves with speed $S_L$, $S_L^\star$, $S_M$, $S_R^\star$ e $S_R$, as 
discussed in the Appendix A. 
The states ${\bf Q}^\star$ and ${\bf Q}^{\star\star}$ are defined as 
\[{\bf Q}^\star_\alpha = (\rho^{\star}_\alpha, E^{\star}_\alpha, \rho^{\star}_\alpha {\bf u}^{\star}_\alpha,{\bf B}^{\star}_\alpha, \psi^{\star}_\alpha) \text{ and } {\bf Q}^{\star\star}_\alpha = (\rho^{\star\star}_\alpha, E^{\star\star}_\alpha, \rho^{\star\star}_\alpha {\bf u}^{\star\star}_\alpha,{\bf B}^{\star\star}_\alpha, \psi^{\star\star}_\alpha),\] 
with $\alpha$ denoting left ($L$) or right ($R$) states. 
\\ \\
\item  The same procedure is done for $B_y$ in the $y$-direction.
\\
\item The computed values of $\psi$ are used to update the mixed correction source term for $\psi^{n+1}$, computing $\psi^{n+1}= \exp \left(-\Delta t^n \frac{c_h^2}{cp^2} \right)\, \psi$.
\end{enumerate}

The adaptive Multiresolution (MR) method of the present paper has been designed to speed
up  finite volume schemes for conservation laws. 
In the following, a brief summary of this technique is given. 
For  a detailed description of these strategies, we refer to \cite{RSTB03,DominguesRousselSchneider:IJNME2009,DominguesGomesRousselSchneider:APNUM2009,Domingues20083758,DGRSESAIM:2011}.

The key ingredient of MR schemes is the decay properties of the wavelet coefficients of the numerical solution. The decay rate indicates the local regularity of the solution. 
In regions where the solution is smooth the coefficients are of small magnitude and thus coarser meshes can be used. 
In regions where the coefficients are significant the numerical solution is less smooth and strong gradients or even jumps are present and a fine mesh must be used \cite{CohenKaberMullerPostel:2003}.
Stopping the refinement in a cell at a certain scale level, where the wavelet coefficients are non-significant leads to an adaptive MR representation.

\begin{figure}[htb]
\begin{center}
\begin{tabular}{c}
\includegraphics[width=0.3\linewidth]{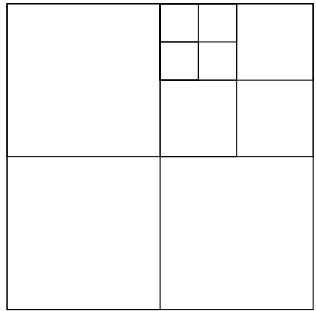}\\[-4mm]
\end{tabular}
\end{center}
\caption{Example of a zoom in a dyadic adaptive Cartesian mesh. Regions where the mesh is refined are associated with detected structures in the solution, \textit{i.e.}, where the wavelet coefficients are significant.
}
\label{fig:mesh}
\end{figure}

For a finite volume scheme the uniform cell-average representation is replaced by cell-averages on an adaptive locally refined mesh, which is formed by the cells whose wavelet coefficients are significant and above a given threshold. An example of an adaptive Cartesian mesh is presented in Fig.~\ref{fig:mesh}.
In MHD solutions localized structures are present, such as discontinuities or shocks. 
They could appear in different space positions in different variables. Thus, the adaptive mesh of the MHD system is a union of the individual adaptive meshes of each quantity.

Tree structures are the natural way to store the reduced MR data.
Mesh adaptivity is then related to an incomplete tree and the refinement can be interrupted at intermediate scale levels.
In other words, using the tree terminology, a MR mesh is formed by leaves, which are nodes without children. 
These leaves correspond to the cell which is being  evolved in time.
In summary, there are three steps in the application of a MR scheme: refinement, evolution, and coarsening.
The refinement operator accounts for possible translations of the solution or the creation of finer scales in the solution between two subsequent time steps. Since the localized structures and thus the local regularity of the solution may change with time, the MR mesh at time $t^n$ may not be sufficient any more at the next time step $t^{n+1}$.
Hence, before evolving the solution in time, the representation of the solution should be interpolated onto an extended mesh that is expected to be a refinement of the adaptive mesh at $t^n$, and to contain
the adaptive mesh at $t^{n+1}$. After that, the time evolution operator is applied to the leaves of the extended mesh.
The numerical fluxes between cells of different levels are computed by adding extra cells, 
called virtual leaves, which will however not be used in the time evolution. 
Conservation is ensured by the fact that the fluxes are always computed on a higher level, the value being projected onto the leaves of a lower level. 
Then, wavelet thresholding is applied in order to unrefine the cells in the extended grid (coarsening) that are not necessary for an accurate representation of the solution at $t^{n+1}$.
This data compression is based on the definition of deletable cells, where the wavelet coefficients which are not significant, \textit{i.e.}, their magnitudes are below a threshold parameter $\epsilon^\ell$, where $\ell$ denotes the cell scale level, are called deletable cells. 
The data compression is the given by
\begin{equation}
D_c=\frac{100\,\sum\limits_{i=1}^{N}C_n(i)}{2^L\,N},\nonumber
\end{equation}
where $N$ is the total number of iterations and $C_n(i)$ is the number of cells in the adaptive mesh at iteration $i\in\{1,\cdots, N\}$.
The number of cells on the finest mesh is defined as $2^L$, where $L$ the finest scale level. 
However, to compute the flux in a conservative form, additional neighbor cells at the same level are also necessary. These neighbor cells are not necessarily present on the adaptive mesh. Thus, if this is the case, we add these neighbor cells to the adaptive mesh, nevertheless they are not evolved in time. Therefore, the memory ised is the sum of the cells of the adaptive mesh plus these neighbor cells. More details in \cite{Roussel:2003,RSTB03}. 

In order to control the $L^1$-norm, Harten's thresholding strategy is used, where
\begin{equation}
\epsilon^\ell=\frac{\epsilon^0}{|\Omega|} 2^{d(\ell-L+1)}, \;\;0\leq \ell \leq L-1,
\label{eq:MR}
\end{equation}
and $d=2$ is the space dimension and, in this two-dimensional case $|\Omega|$ is the area of the domain.
Therefore, in the Harten's strategy, we use a smaller value of the parameter $\epsilon$ in the coarser scales than in fines scales.  
For comparison, we shall also consider level independent threshold parameters: $\epsilon^\ell = \epsilon$, for all $\ell$. 
Herein, the multiresolution analysis corresponds to a prediction operator based on a third order polynomial 
interpolation on the cell-averages \cite{RSTB03}. 
We recall that time integration is performed by a second order Runge--Kutta scheme.

\section{Numerical experiments}
\label{sec:numerical}


We present here a 2D Riemann numerical experiment to illustrate the efficacy of our method compared to the traditional FV scheme. 
For the  2D Riemann initial condition we have used the values of the MHD variables presented in Table~\ref{table:R2D}.
The computational domain is $[-1,1]\times[-1,1]$ and  Neumann boundary conditions have been applied. 
This example is proposed in \cite{Dedneretal:2002}, except for the boundary condition.

We have also chosen  $\gamma=5/3$, the final time of computations $t=0.1$ and $t=0.25$, the CFL parameter $C_{CFL}=0.3$ and $c_p^2/c_h=0.18$. 
We have tested $\epsilon^\ell=\epsilon=0.010,0.008,0.005$ and Equation~\ref{eq:MR} with $\epsilon^0=0.05,0.03,0.01$.

The reference GLM--MHD FV code used in this work has been developed in $C^{++}$ language, inspired by the Fortran code developed by \cite{Bastien:2009}, including an upgrade and new features for the implementation of the numerical flux HLLD.
The GLM--MHD MR code developed  in \cite{Gomes:2012:AnMuAd} is based on the hydrodynamics MR Carmen code developed in \cite{RSTB03, Roussel:2003}.The implementation has been optimized improving the momory allocation and unrolling the for-loops for the allocation of the variables. The CPU is improved about a factor 4for the test case studied here with $L=8$ adaptive scales and $\epsilon^0=0.01$.

For the numerical error analysis we have used a reference solution computed with a GLM--MHD FV scheme with $L=11$ scales using the same numerical scheme in space, implemented in the AMROC code \cite{Deiterdingetal:2009} which is parallelized. We computed the $L_1$-error for the density solution ($L_1^e(\rho)$). The CPU time for the MHD-FV reference is obtained with another code that is not parallel.

\begin{table}[htb]
  \caption{Initial condition of the 2D Riemann problem. The domain is $[-1,1]\times[-1,1]$ with Neumann boundary conditions and $\gamma=\dfrac{5}{3}$.}
  \label{table:R2D}
  \begin{center}
 \begin{small}
  \begin{tabular}{ccccp{0.1mm}cccc}
    \hline   \multicolumn{9}{c}{$\boldsymbol{x>0}$} \\
    \hline
      \multicolumn{4}{c}{$\boldsymbol{y<0}$} && \multicolumn{4}{c}{$\boldsymbol{y>0}$}\\
    \cline{1-4}\cline{6-9}     
     $\rho$ & $\rho\,u_x$ & $\rho\,u_y$ & $\rho\,u_z$ && $\rho$ & $\rho\,u_x$ & $\rho\,u_y$ & $\rho\,u_z$\\
    1.0304 & 1.5774 & -1.0455 & -0.1016& & 0.9308 & 1.4557 & -0.4633  & 0.0575\\
    $E$ & $B_x$ & $B_y$ & $B_z$ && $E$ & $B_x$ & $B_y$ & $B_z$\\
    5.7813 & 0.3501 & 0.5078 & 0.1576 && 5.0838 & 0.3501 & 0.9830 & 0.3050\\[1mm]
    \hline   \multicolumn{9}{c}{$\boldsymbol{x<0}$} \\
    \hline
       \multicolumn{4}{c}{$\boldsymbol{y<0}$} && \multicolumn{4}{c}{$\boldsymbol{y>0}$}\\
    \cline{1-4}\cline{6-9}    
     $\rho$ & $\rho\,u_x$ & $\rho\,u_y$ & $\rho\,u_z$ && $\rho$ & $\rho\,u_x$ & $\rho\,u_y$ & $\rho\,u_z$\\
       1.0000 & 1.7500 & -1.0000 & 0.0000 &&1.8887 & 0.2334 & -1.7422 & 0.0733\\
     $E$ & $B_x$ & $B_y$ & $B_z$ && $E$ & $B_x$ & $B_y$ & $B_z$\\
        6.0000 & 0.5642 & 0.5078 & 0.2539 && 12.999 & 0.5642 & 0.9830 & 0.4915\\
	\hline
  \end{tabular}
  \end{small}
\end{center}
   \end{table}

\begin{figure}[htb]
\begin{center}
\begin{tabular}{cc}
$\rho$ & $B_y$ \\
\includegraphics[width=0.45\linewidth]{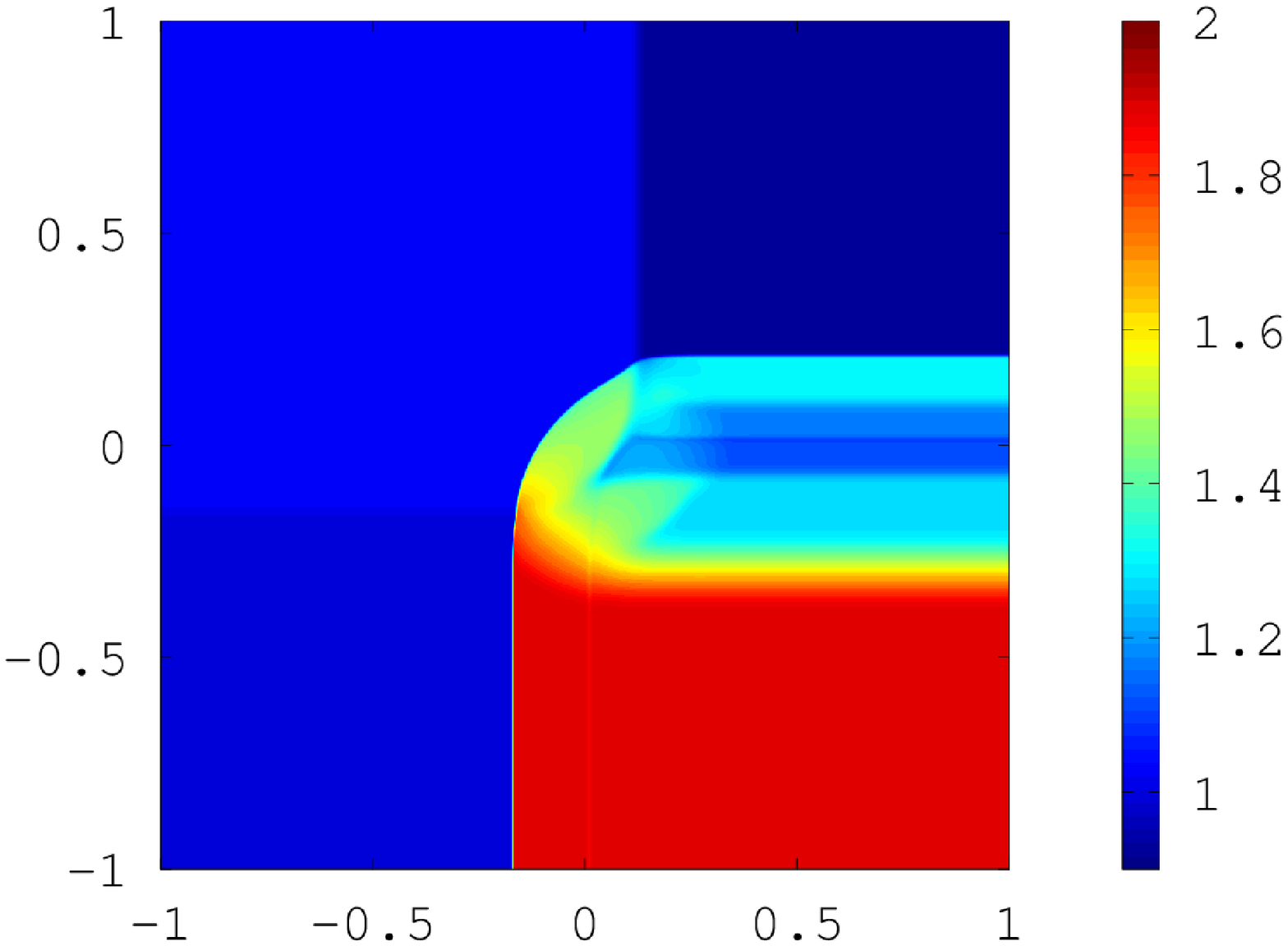} & \includegraphics[width=0.45\linewidth]{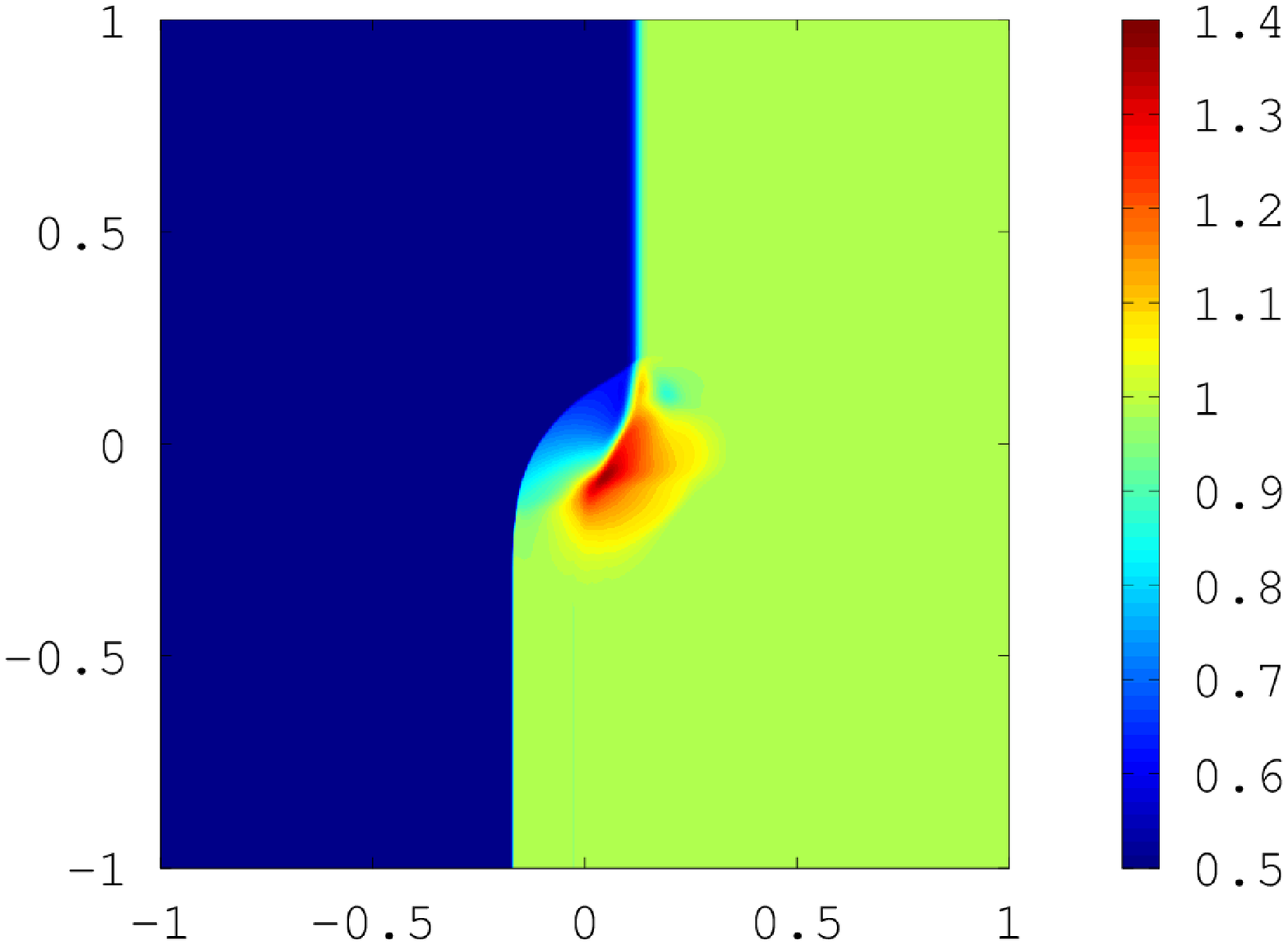}\\
$u_y$ & $u_z$ \\
\includegraphics[width=0.45\linewidth]{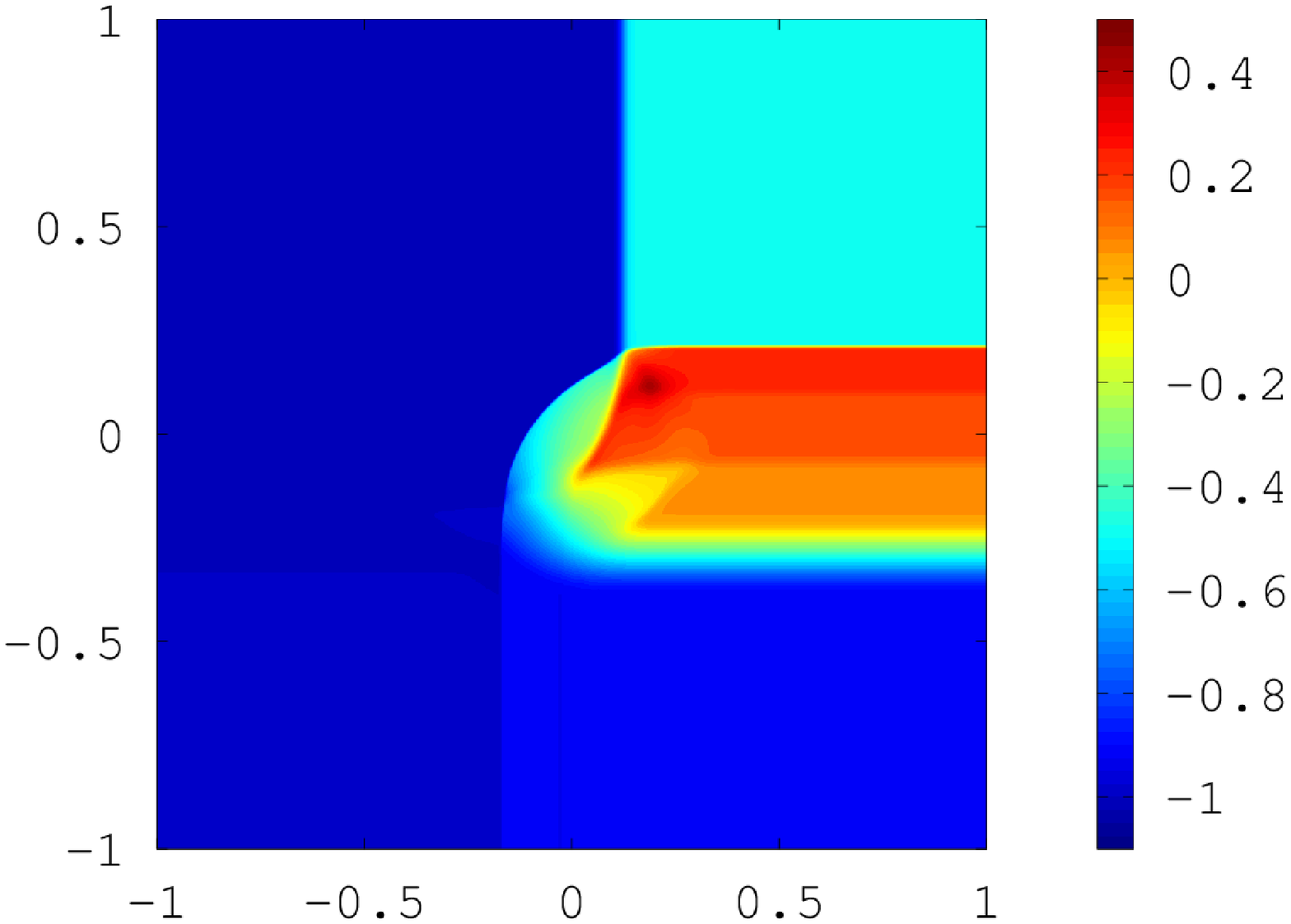} & \includegraphics[width=0.45\linewidth]{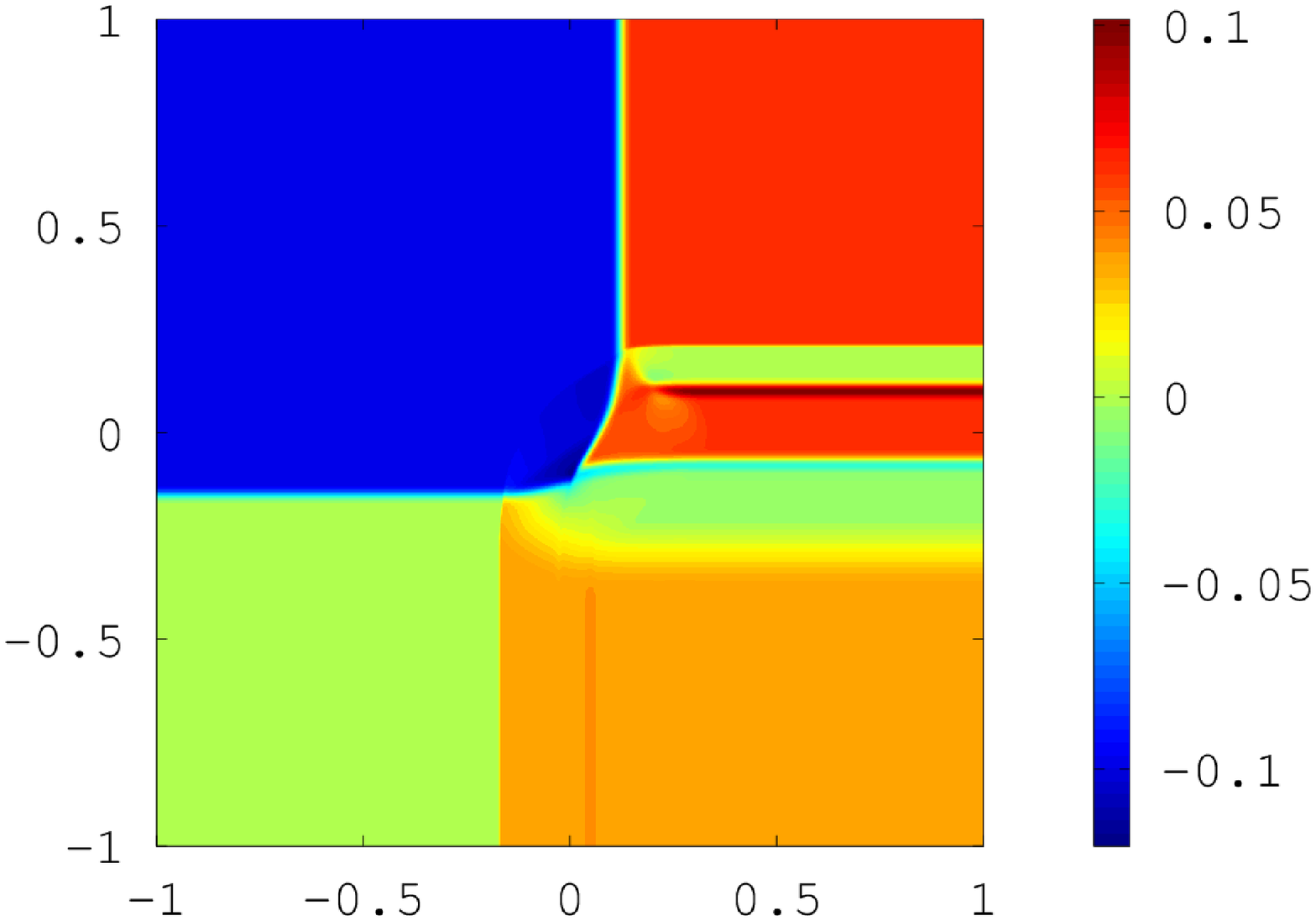}
\end{tabular}
\end{center}
\caption{FV reference solution for the 2D Riemann problem using GLM--MHD with mixed correction. Shown are variables $\rho$, $B_y$, $u_y$ and $u_z$ obtained at time $t=0.1$ and $L=11$. 
}
\label{fig:2DRref}
\end{figure}

\begin{figure}[htb]
\begin{center}
\begin{tabular}{cc}
$\rho$ & $B_y$ \\
\includegraphics[width=0.45\linewidth]{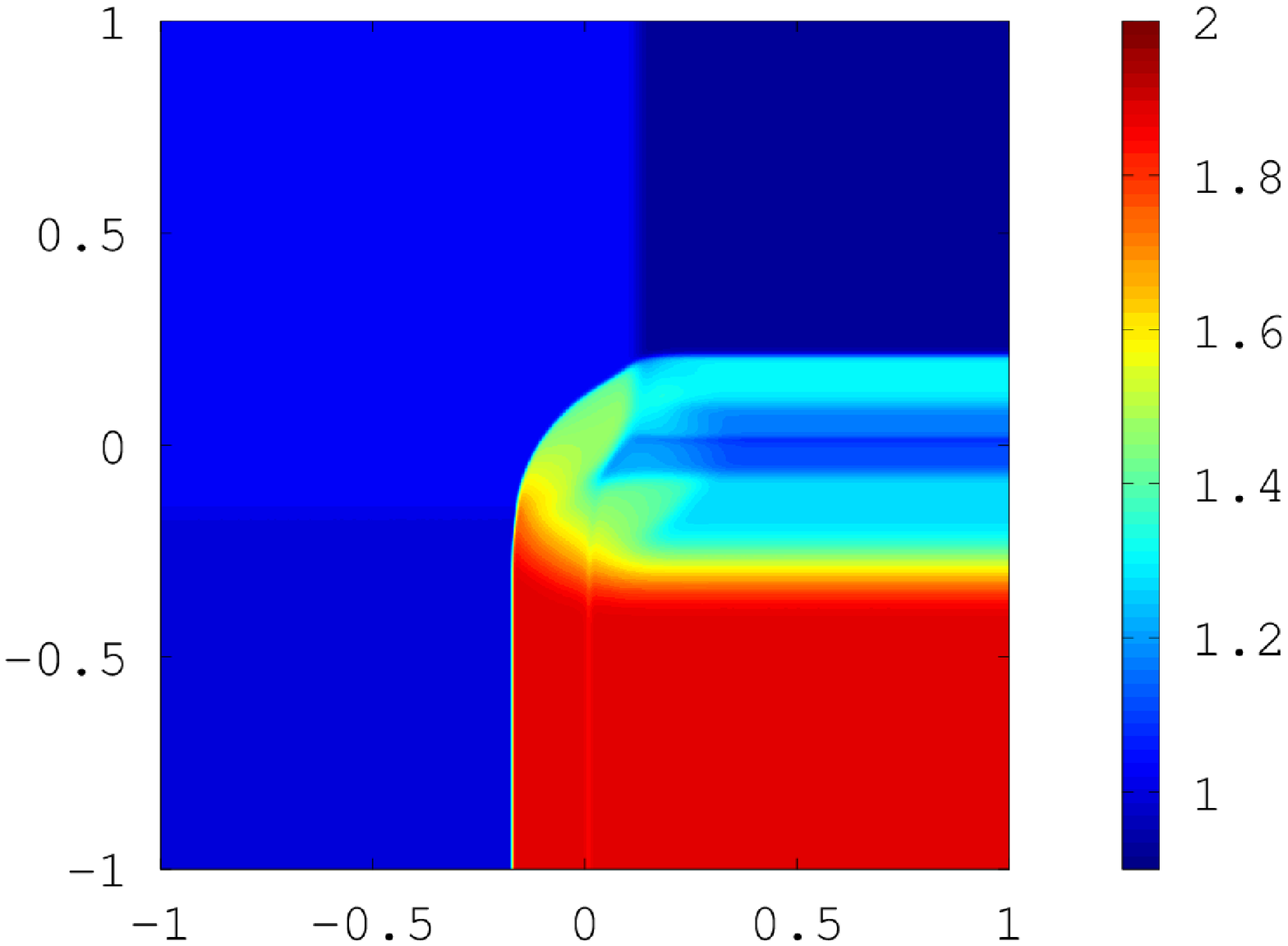} & \includegraphics[width=0.45\linewidth]{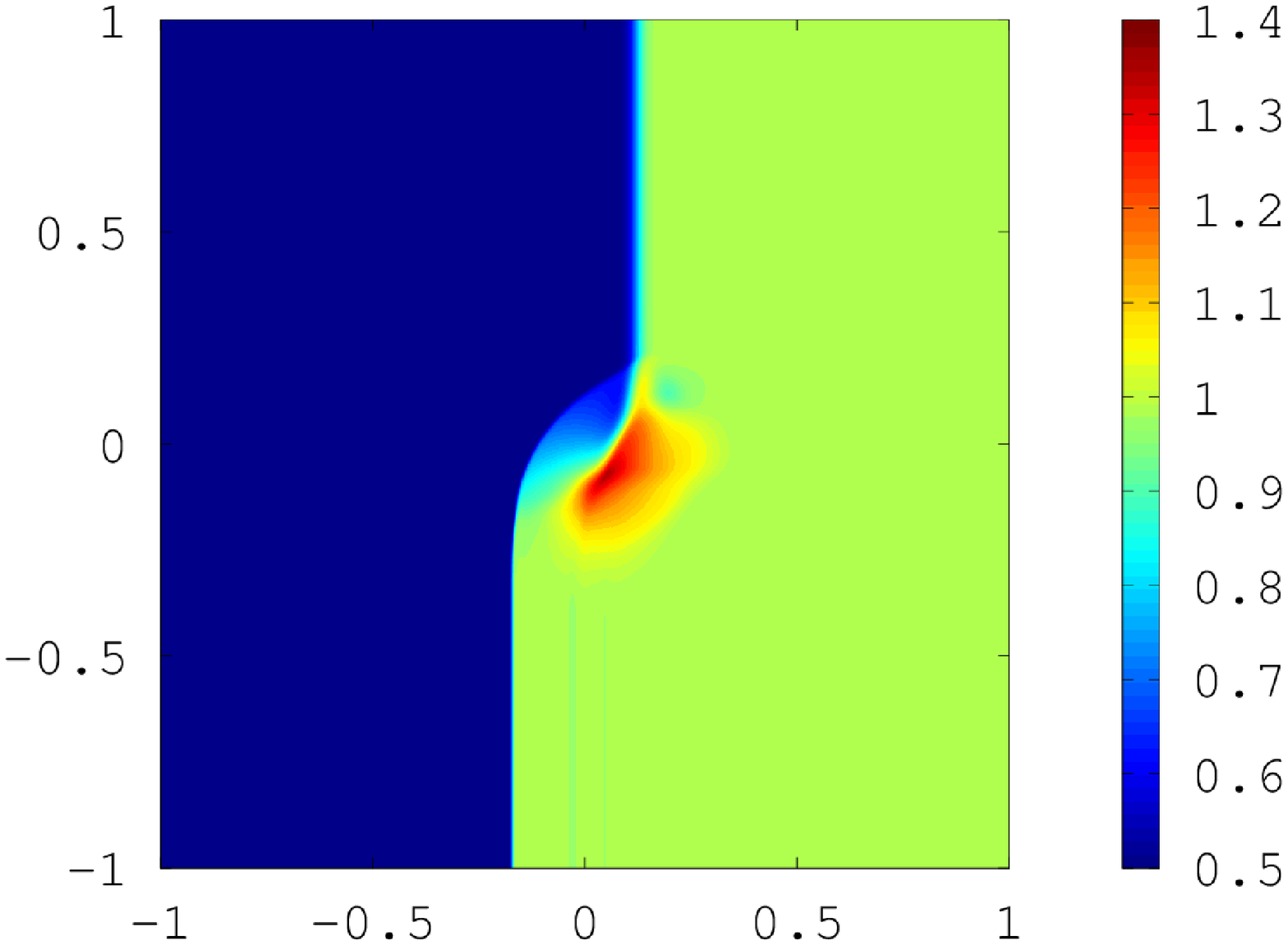}\\
$u_y$ & $u_z$ \\
\includegraphics[width=0.45\linewidth]{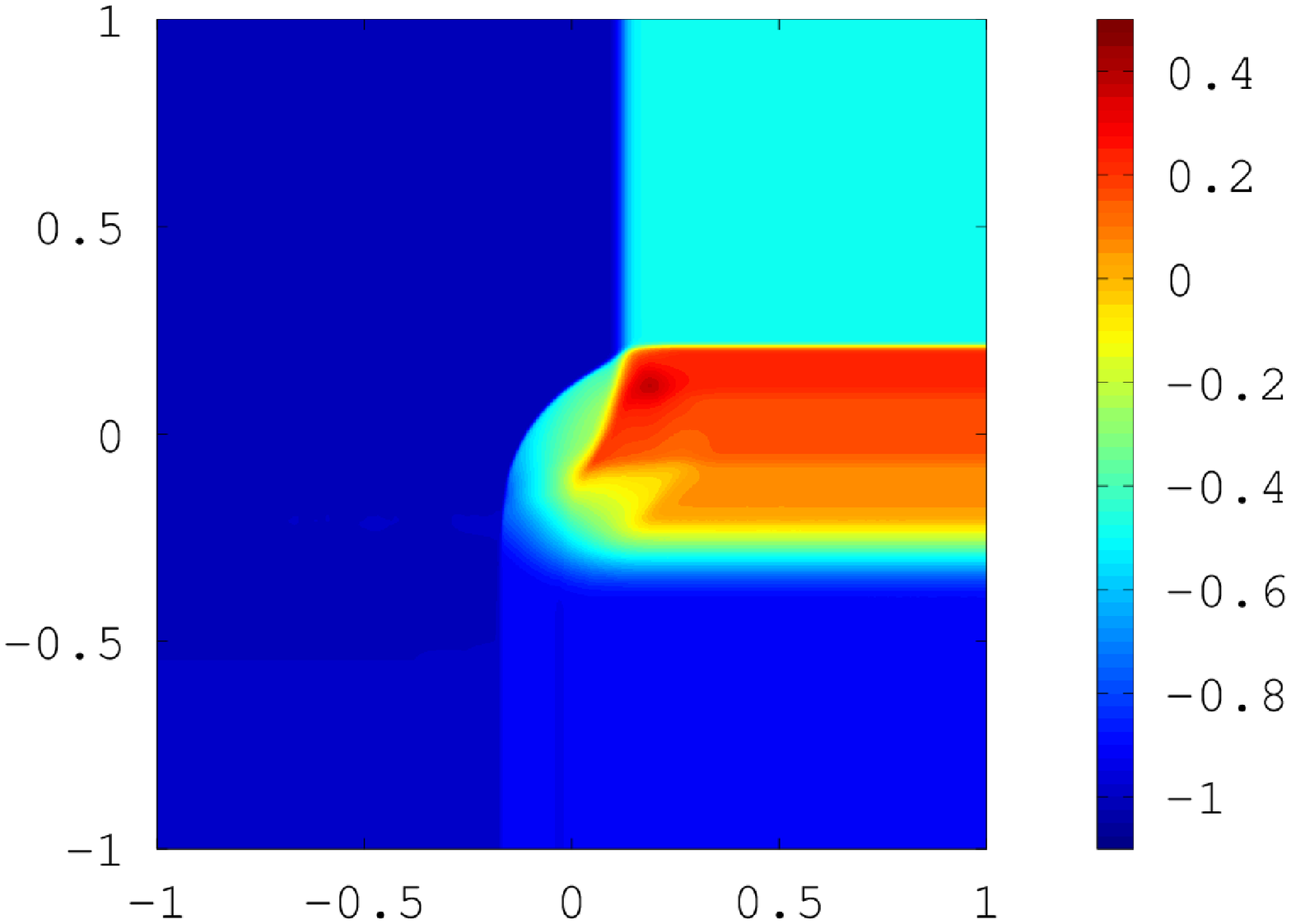} & \includegraphics[width=0.45\linewidth]{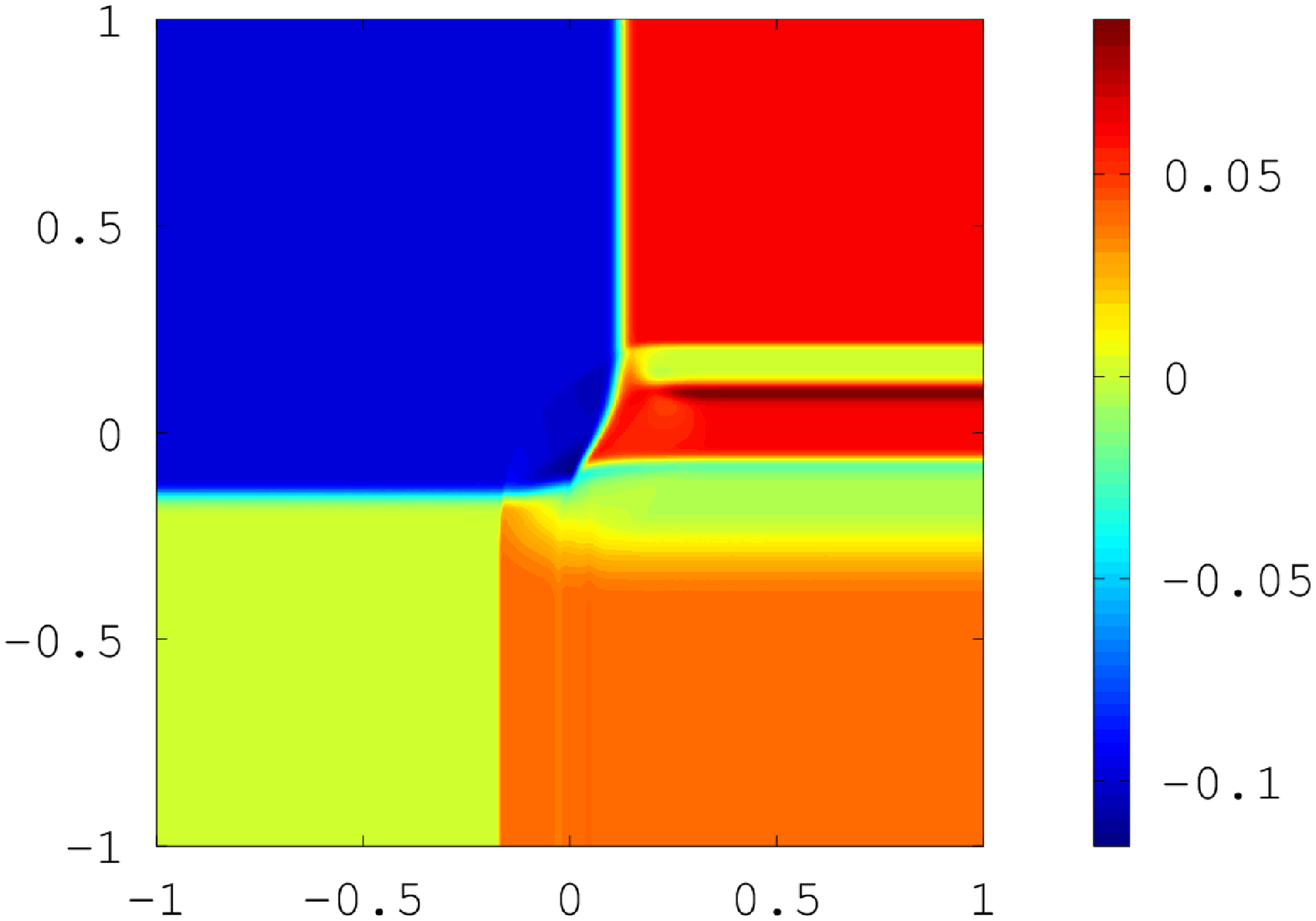}
\end{tabular}
\end{center}
\caption{MR solutions with $\epsilon^0=0.01$ for the 2D Riemann problem using GLM--MHD with mixed correction. Shown are variables $\rho$, $B_y$, $u_y$ and $u_z$ obtained at time $t=0.1$ and $L=10$. 
}
\label{fig:2DRsol}
\end{figure}

\begin{figure}[htb]
\begin{center}
\begin{tabular}{ccc}
$t=0$ & $t=0.1$ & $t=0.25$ \\[-1cm]
\includegraphics[width=0.3\linewidth]{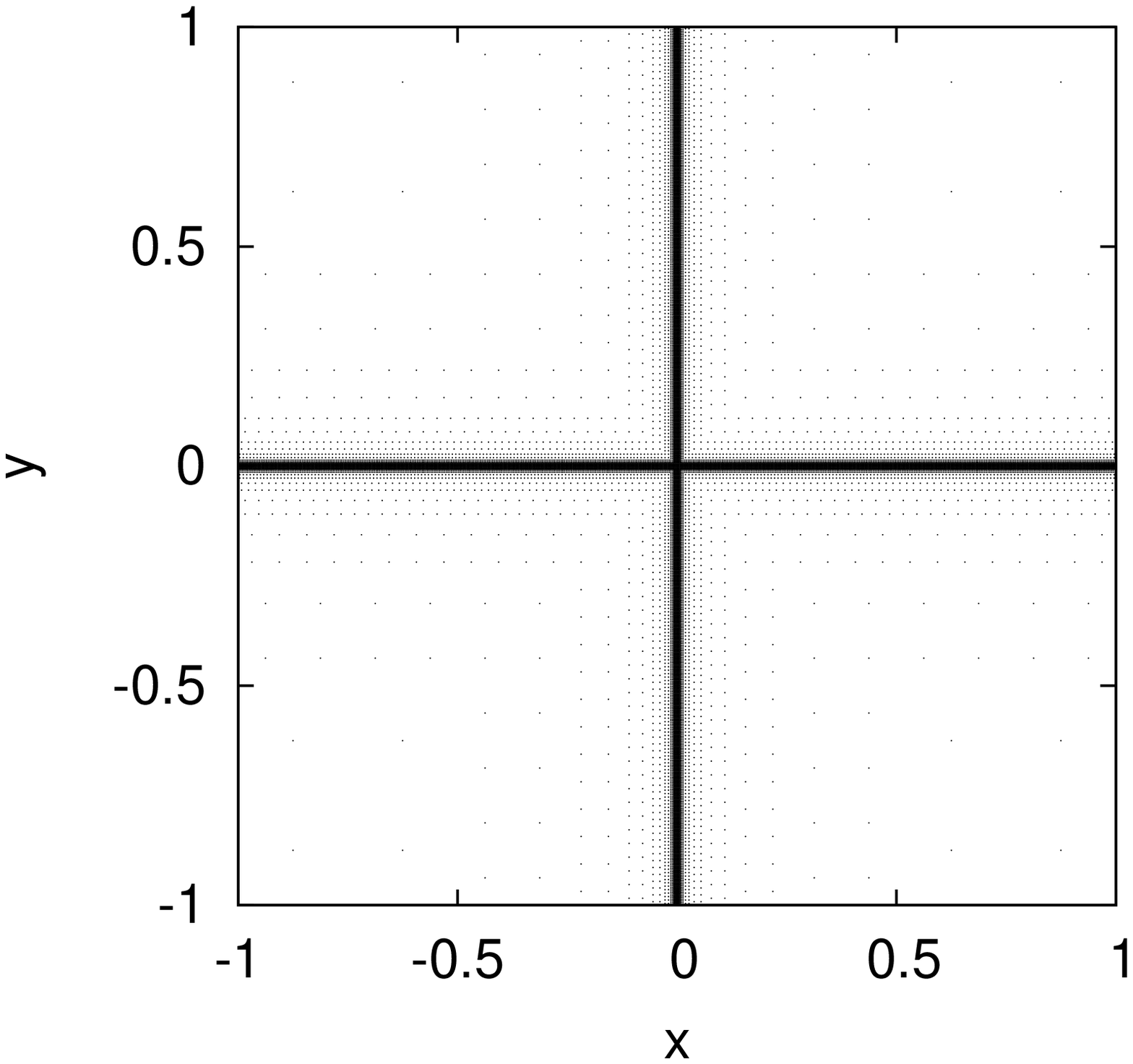} &
\includegraphics[width=0.3\linewidth]{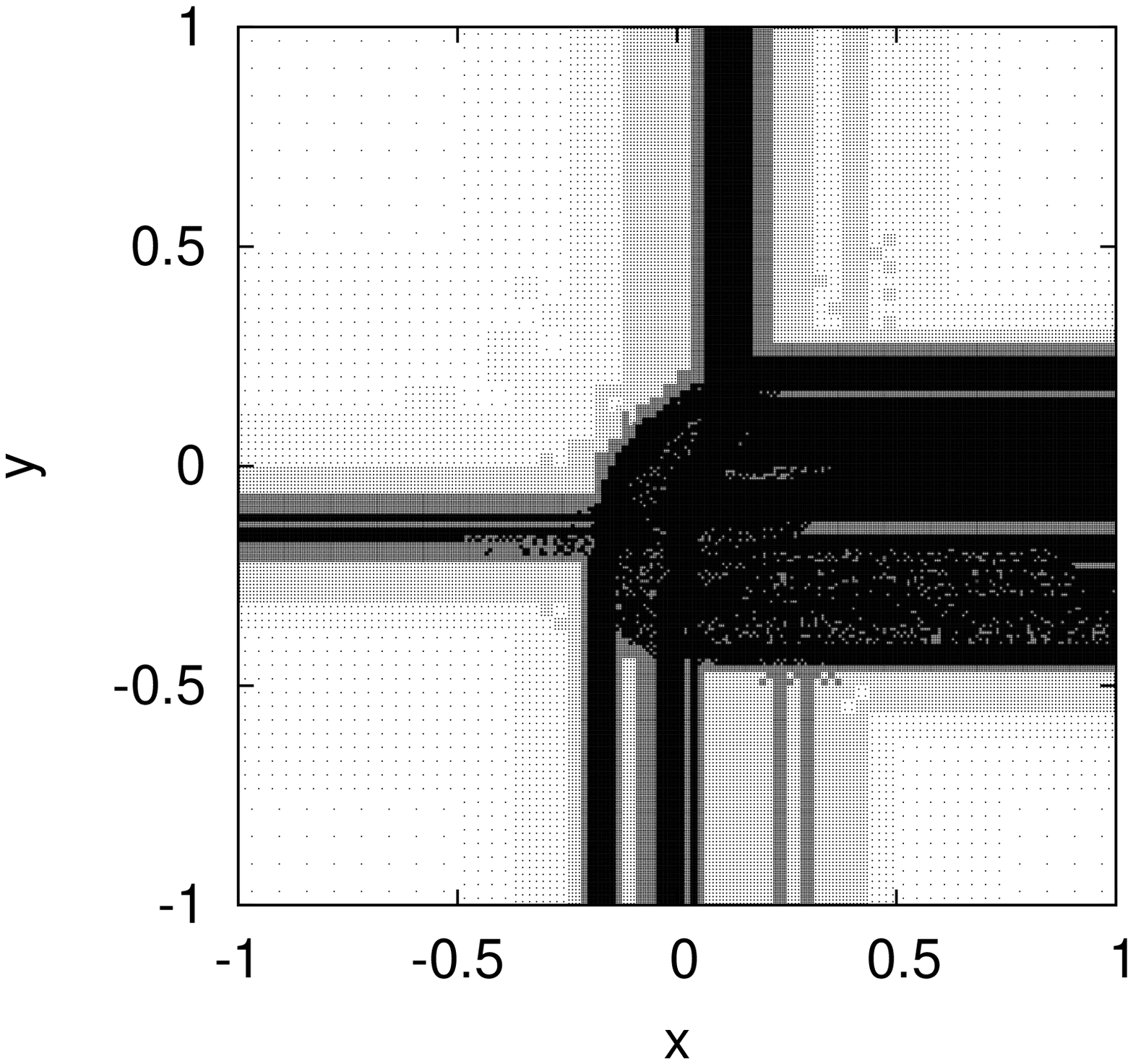} &
\includegraphics[width=0.3\linewidth]{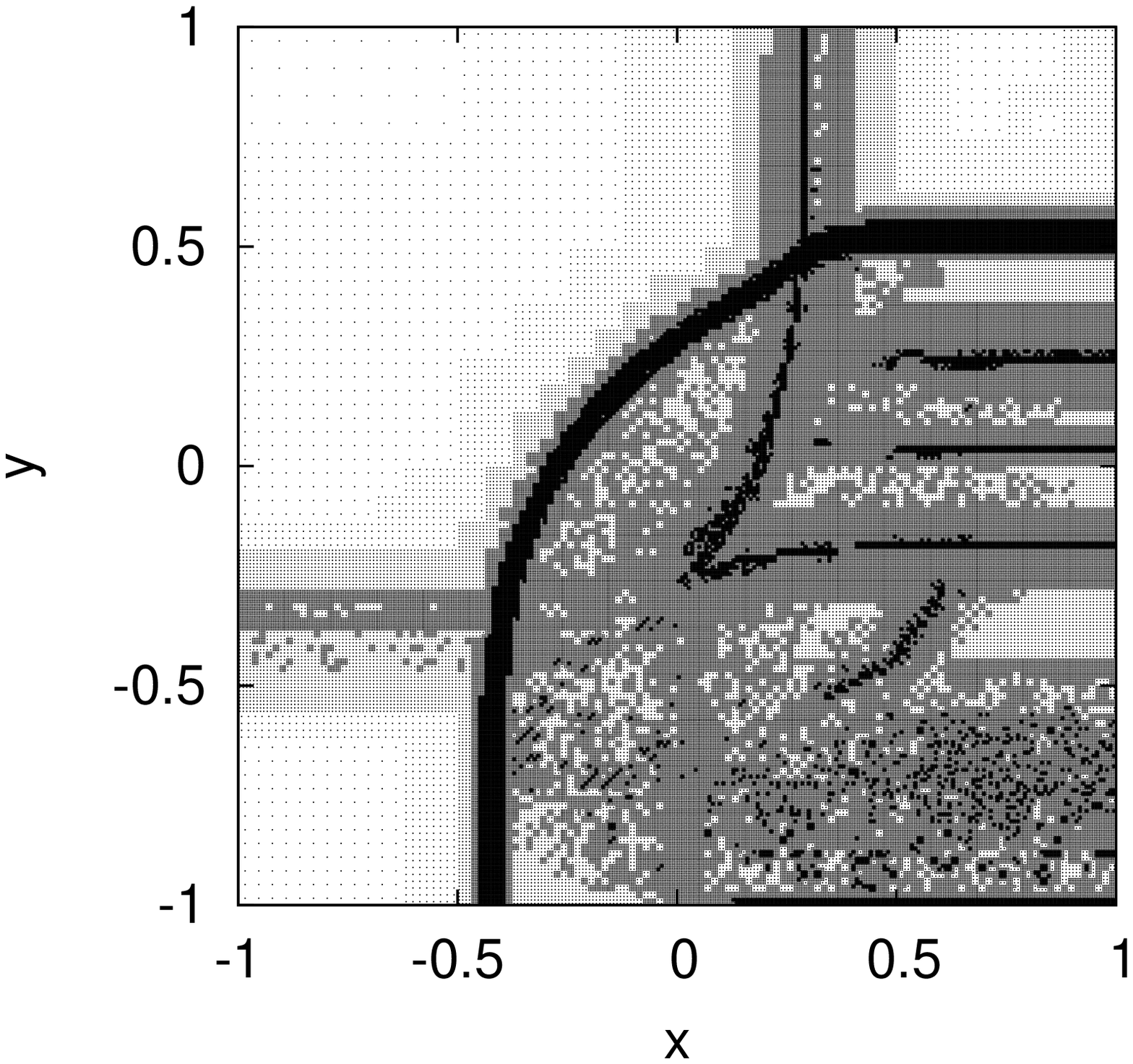}\\[-2cm]
\end{tabular}
\end{center}
\caption{Cell midpoints of the adaptive mesh $L=10$ of the MR computation for the 2D Riemann problem using GLM--MHD with mixed correction at time $t=0$ with $2.30\%$ of the cells, at $t=0.1$ with $26.65\%$ and $\epsilon^0=0.01$; and at time $t=0.25$ with $18.37\%$ of cell and $\epsilon^0=0.05$.
}
\label{fig:2DRmesh}
\end{figure}

\begin{figure}[htb]
\begin{center}
\begin{tabular}{cc}
$\rho$ & $B_y$ \\
\includegraphics[width=0.45\linewidth]{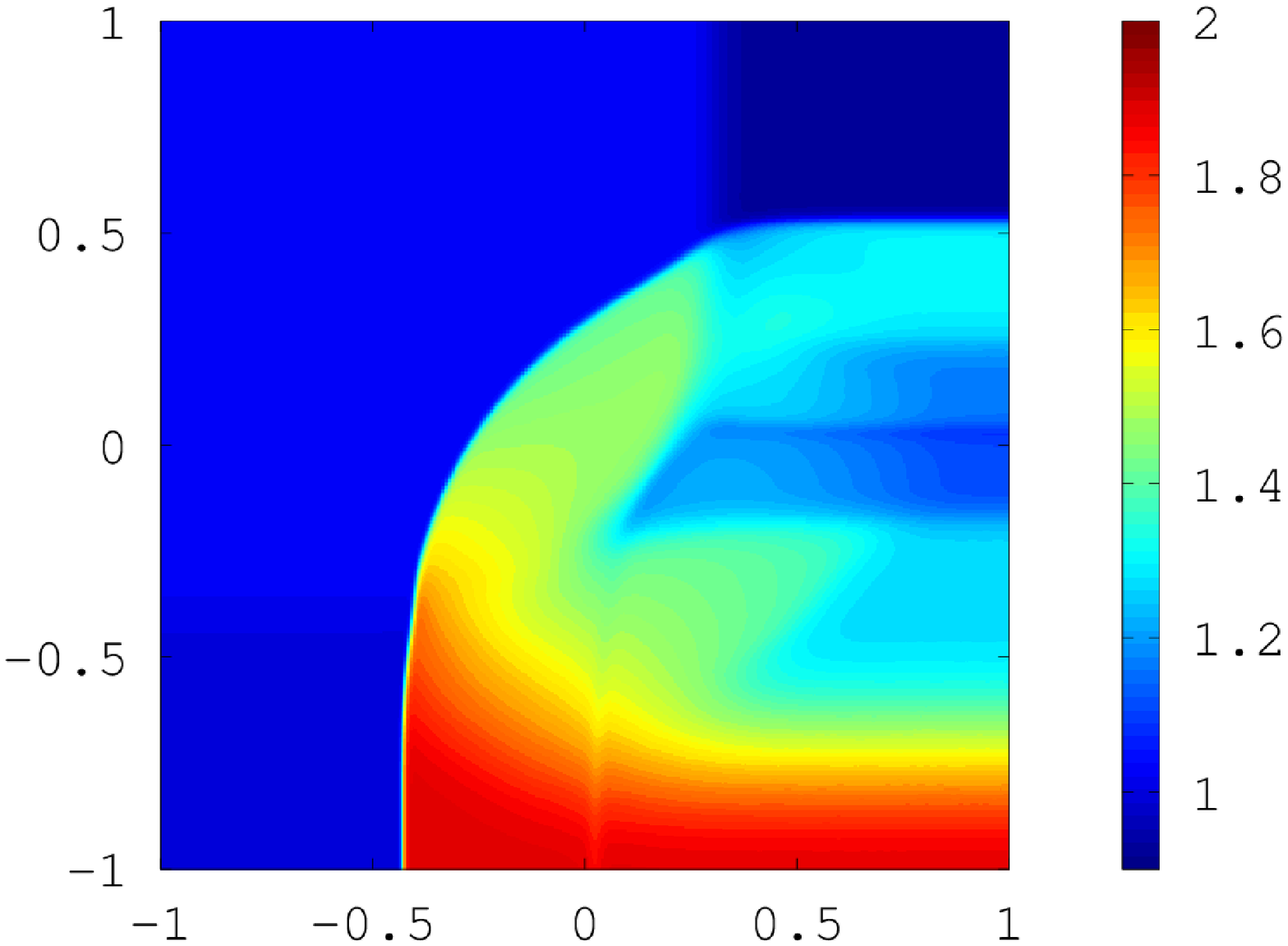} & \includegraphics[width=0.45\linewidth]{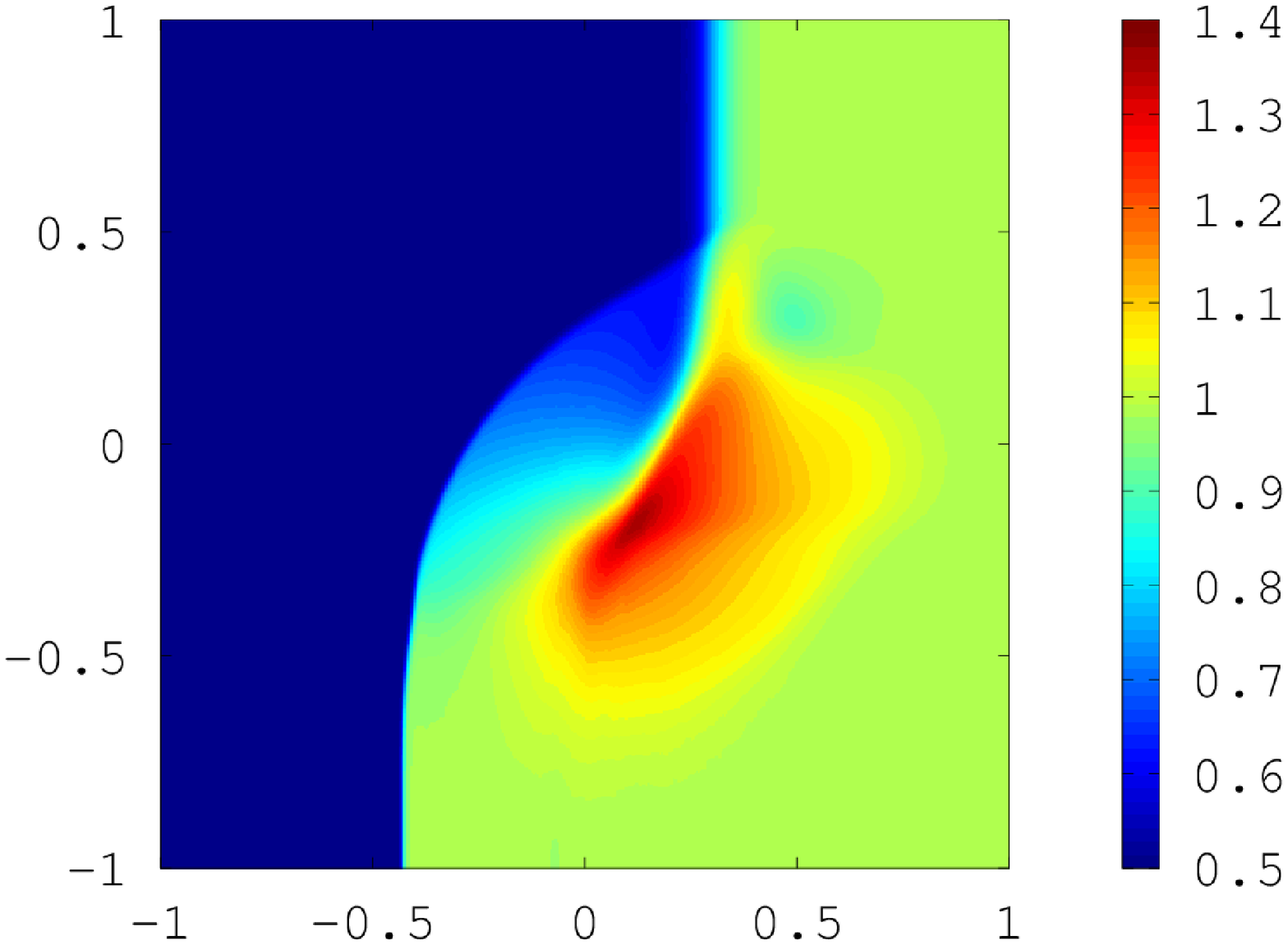}\\
$u_y$ & $u_z$ \\
\includegraphics[width=0.45\linewidth]{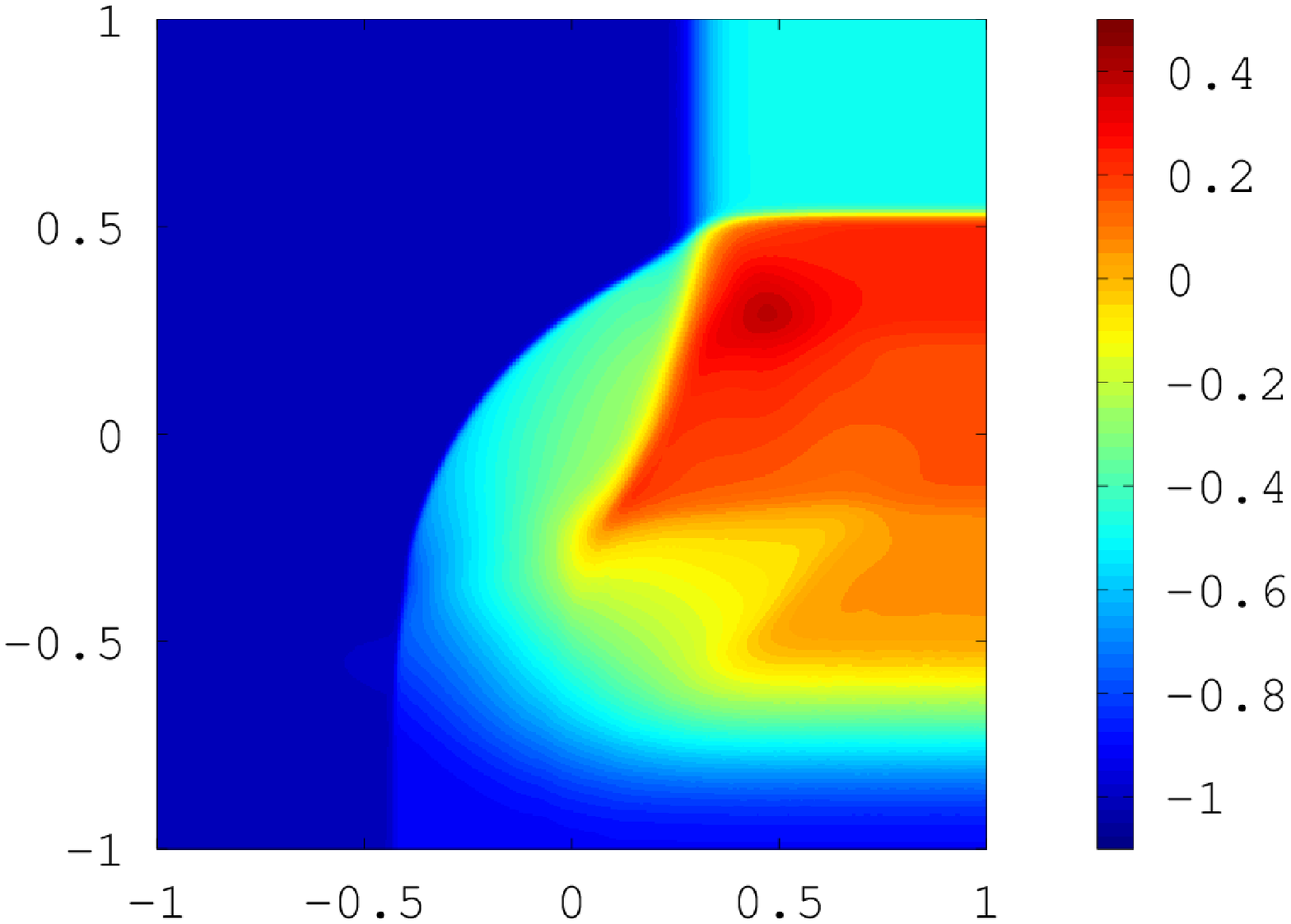} & \includegraphics[width=0.45\linewidth]{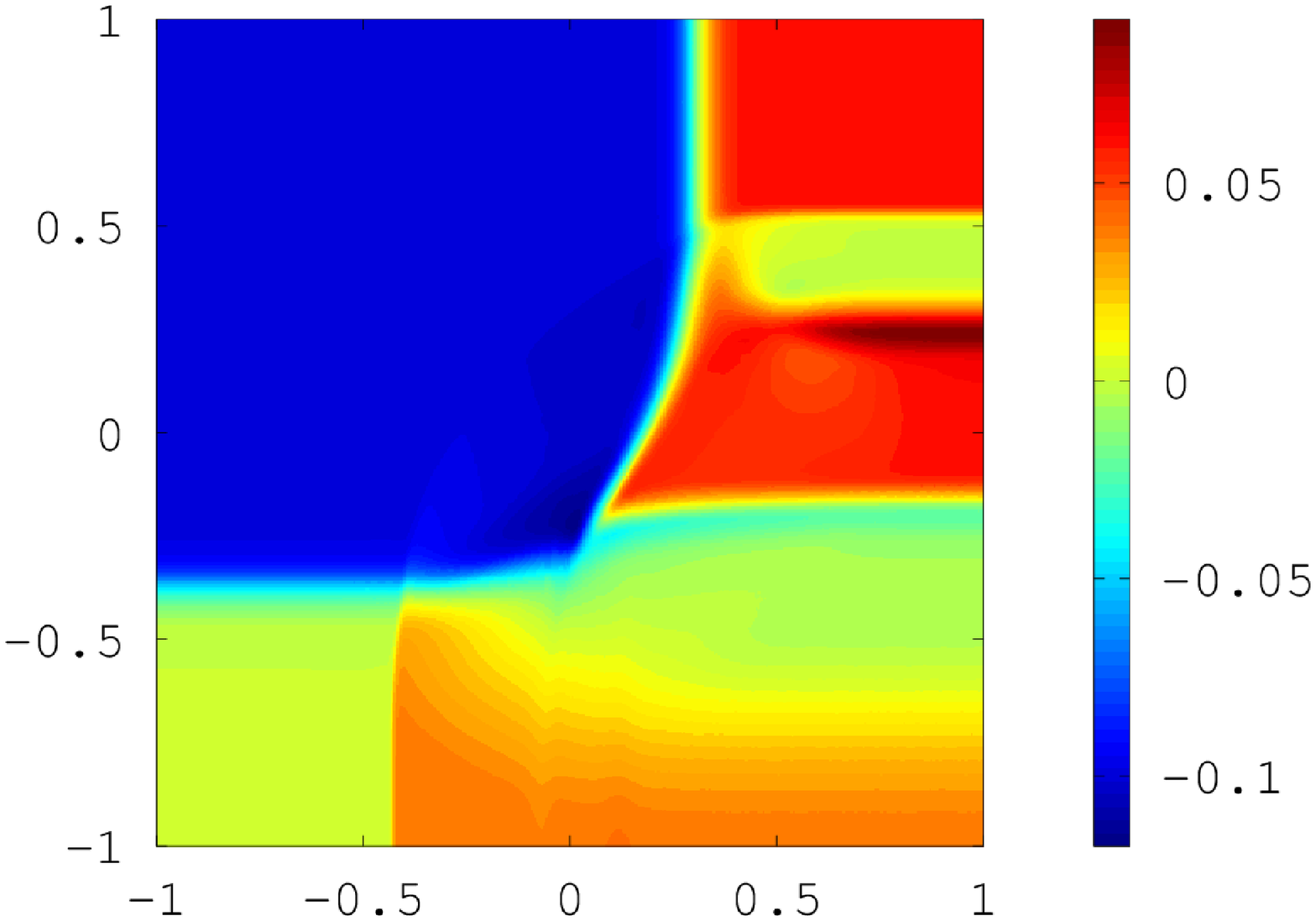}
\end{tabular}
\end{center}
\caption{MR solution for the 2D Riemann problem using GLM--MHD with mixed correction for $\epsilon^0=0.05$. Shown are variables $\rho$, $B_y$, $u_y$ and $u_z$ obtained at time $t=0.25$ and $L=9$. 
}
\label{fig:2DR-solt025}
\end{figure}

The reference solution and numerical MR solutions for $\epsilon^0=0.01$ and $L=10$ at $t=0.1$ are presented in Figs.~\ref{fig:2DRref} and \ref{fig:2DRsol}, respectively. For a later time $t=0.25$, the numerical MR solution with $L=9$ is presented in Fig.~\ref{fig:2DR-solt025}. In the solutions, we can observe that the structures are not always aligned, \textit{e.g.}, we can see a structure that appears in the density but not in the $y$-component of magnetic field in the right part of the domain. In this region, the latter variable is almost constant.
This is expected because in plasma processes the discontinuities may not necessarily occur at the same position for all quantities. 
The $B_x$ component and $p$ (not shown here) have a similar behavior as $\rho$, 
and the $u_z$ component has a similar behavior as $B_z$.  
These observations are expected and they increase the number of cells in the adaptive mesh in the MHD case compared to hydrodynamic case.
Fig.~\ref{fig:2DRmesh} presents an example of the adaptive mesh with $\epsilon^0=0.01$ for the initial, intermediate and final computational time.
We can observe that the adaptive meshes represent  all the structures present in the solutions.

\begin{figure}[htb]
\begin{center}
\begin{tabular}{ccc}
(a) & (b) & (c)\\
\includegraphics[width=0.32\linewidth]{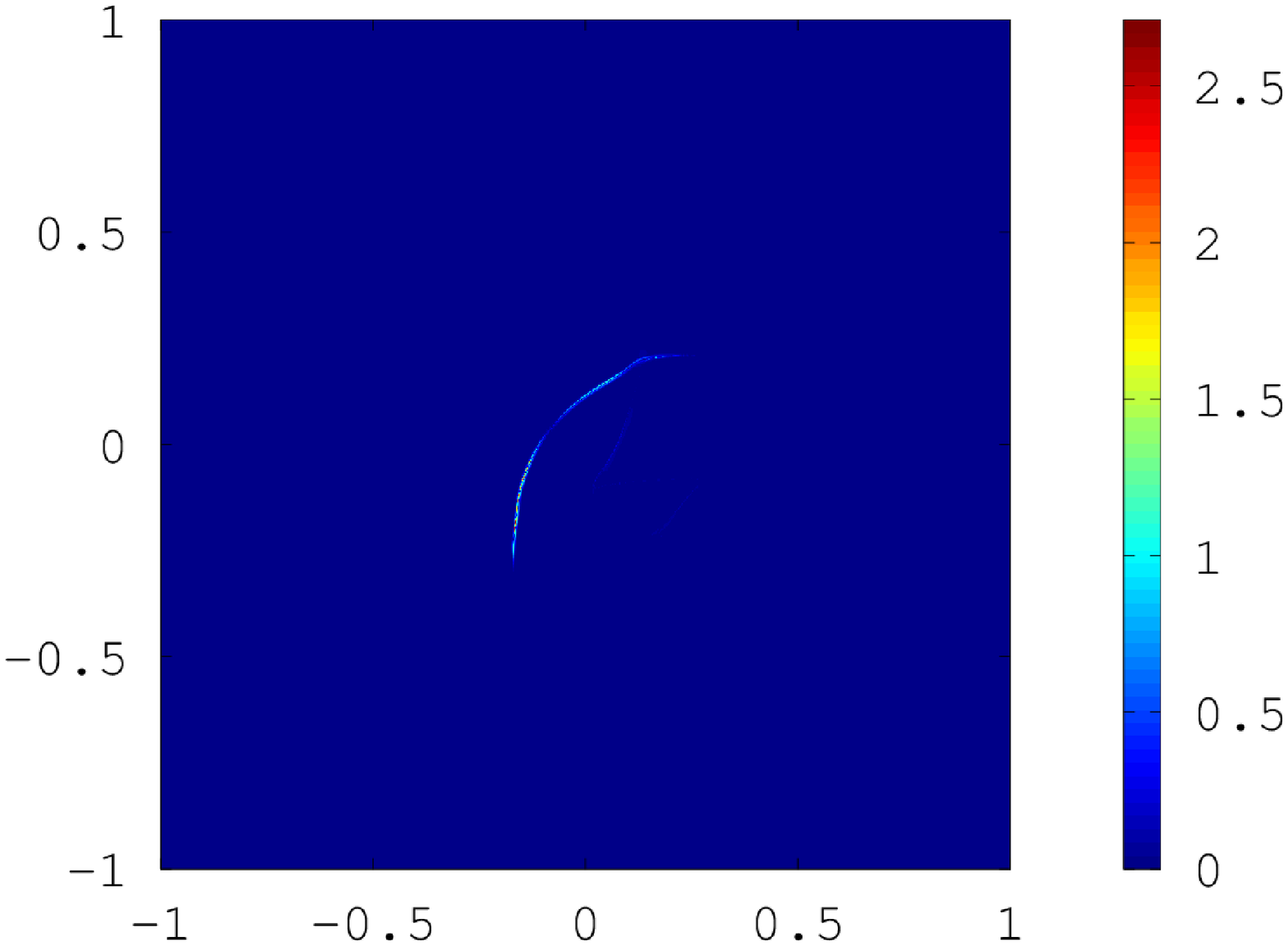} & \includegraphics[width=0.32\linewidth]{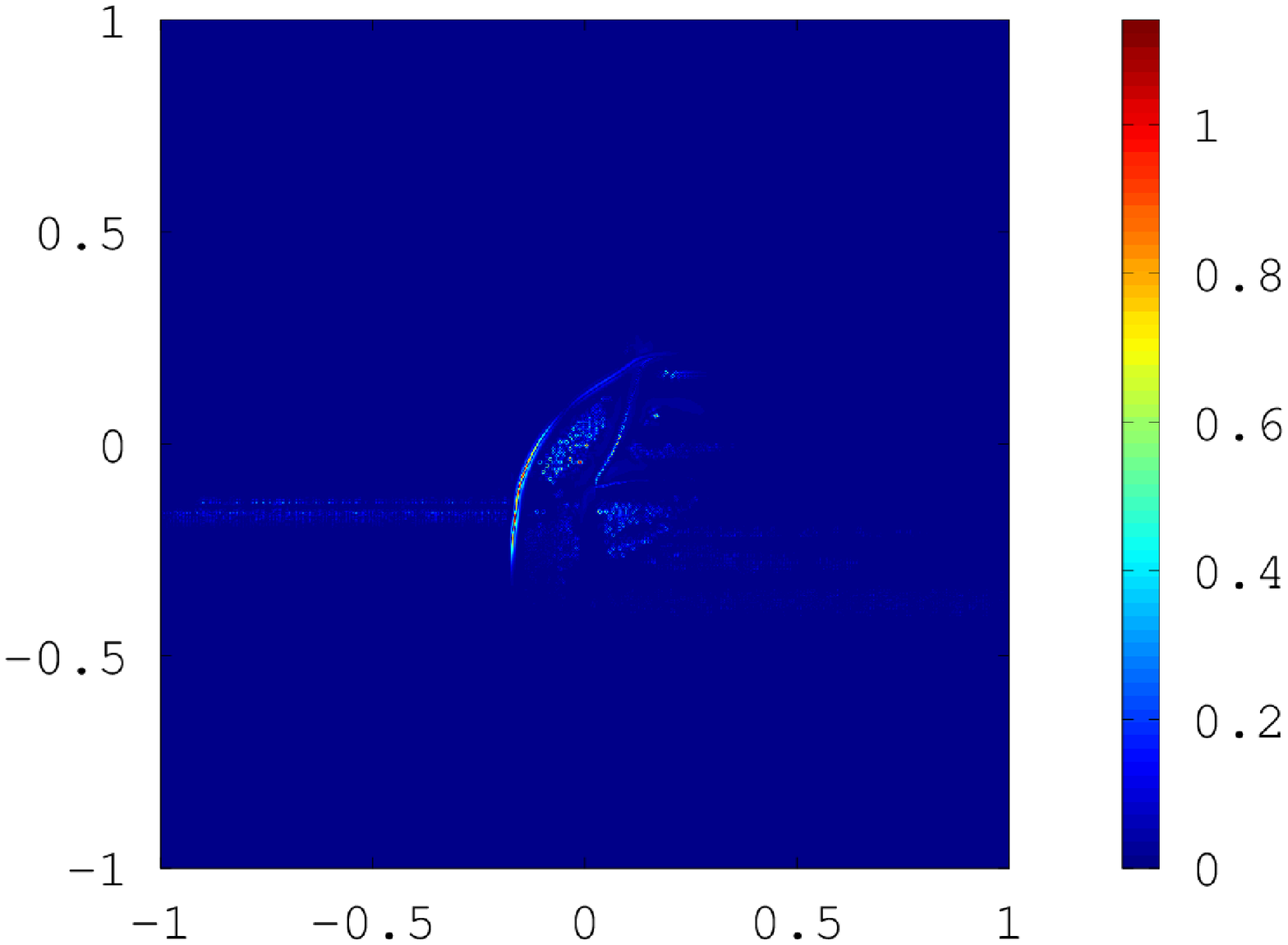} & \includegraphics[width=0.32\linewidth]{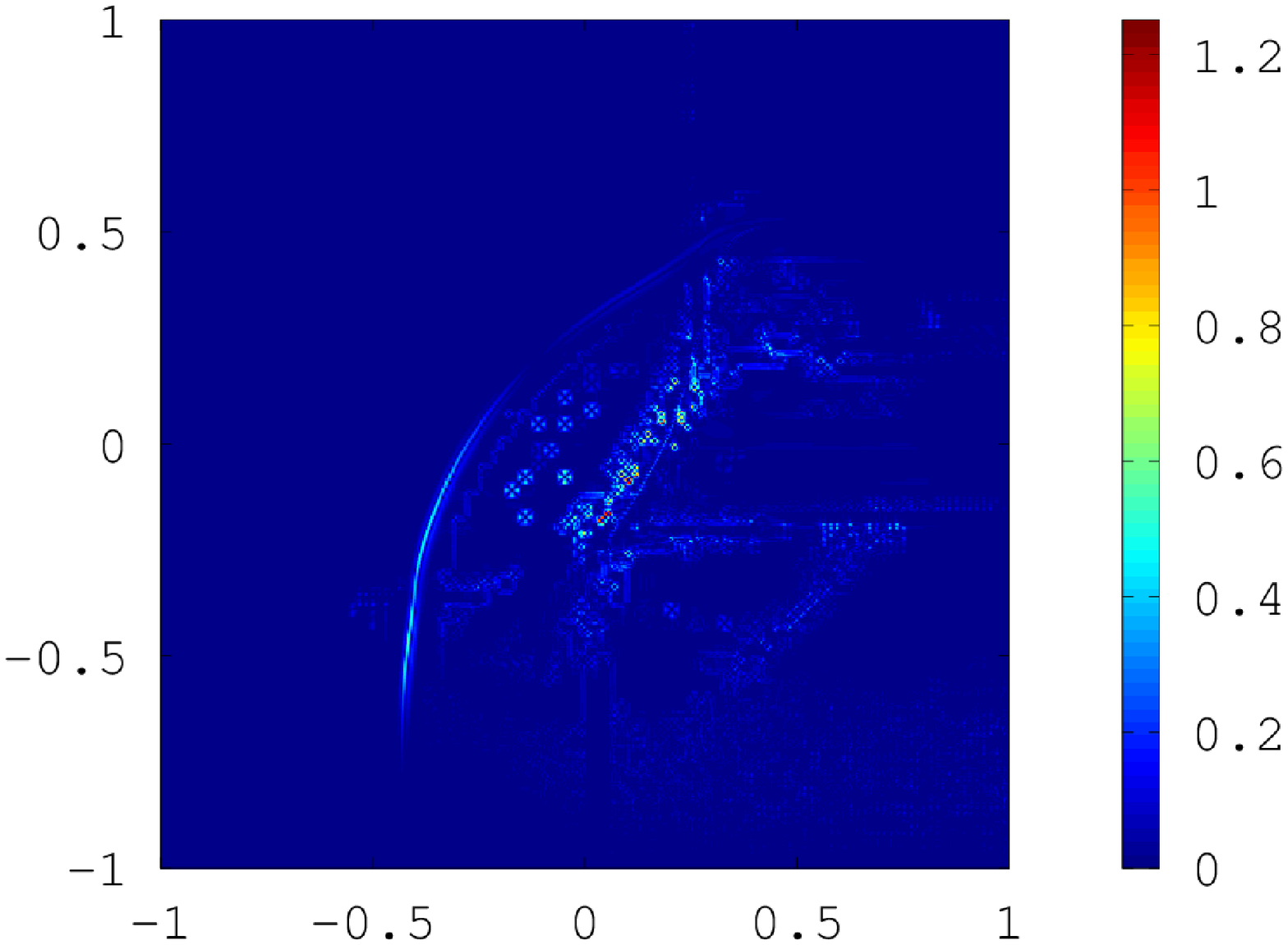}
\end{tabular}
\end{center}
\caption{Values of $\nabla \cdot \bf B$ for the  2D Riemann problem obtained with: $(a)$ FV reference scheme using GLM--MHD with mixed correction for $L=11$; and (b) MR scheme with $\epsilon^0=0.01$ using GLM--MHD with mixed correction for $L=10$ at time $t=0.1$; and (c) MR scheme with $\epsilon^0 =0.05$ using GLM-MHD with mixed correction for $L=10$ at time $t=0.25$. Note that the values of this quantity are mesh-dependent. 
}
\label{fig:divB2DR_image}
\end{figure}

Using the GLM--MHD with the mixed correction, the divergence of the magnetic field is not necessarily zero. 
However, this correction improves the convergence of the numerical solution of the MHD system to the expected physical solution,
as discussed in \cite{Dedneretal:2002}. 
Fig.~\ref{fig:divB2DR_image} presents $\nabla \cdot \bf B$ for the FV reference  for $L=11$ and two MR solutions for $L=10$ with $\epsilon^0=0.01$ at time $t=0.1$ and  $\epsilon^0=0.05$ at time $t=0.25$.
We observe that the maximum values of divergence are in the front transition regions, near the central part of the domain.

\begin{figure}[htb]
\psfrag{epsilon 0}{\tiny{$\epsilon=0$}}
\psfrag{epsilon 0.010}{\tiny{$\epsilon=0.010$}}
\psfrag{epsilon 0.005}{\tiny{$\epsilon=0.005$}}
\psfrag{epsilon 0.008}{\tiny{$\epsilon=0.008$}}
\psfrag{epsilonell 0.010}{\tiny{$\epsilon^0=0.010$}}
\psfrag{epsilonell 0.030}{\tiny{$\epsilon^0=0.030$}}
\psfrag{epsilonell 0.050}{\tiny{$\epsilon^0=0.050$}}
\psfrag{DivBmax}{$B_\mathrm{div}(t)$}
\begin{center}
\begin{tabular}{ccc}
 $(a)\,$MR, $\, L=8$ &  $(b)\,$MR, $\,L=9$\\
\includegraphics[width=0.35\linewidth]{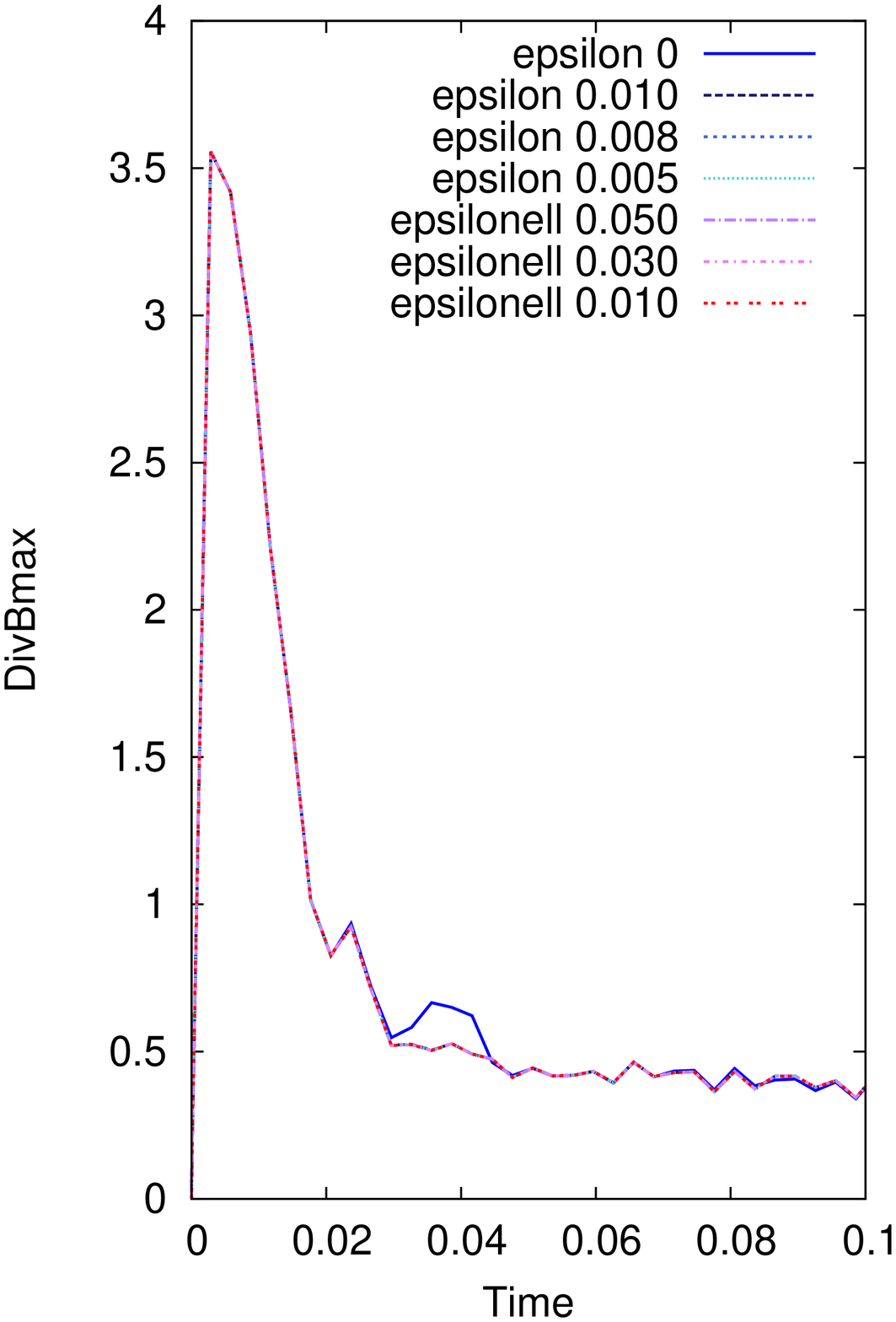}&
\includegraphics[width=0.35\linewidth]{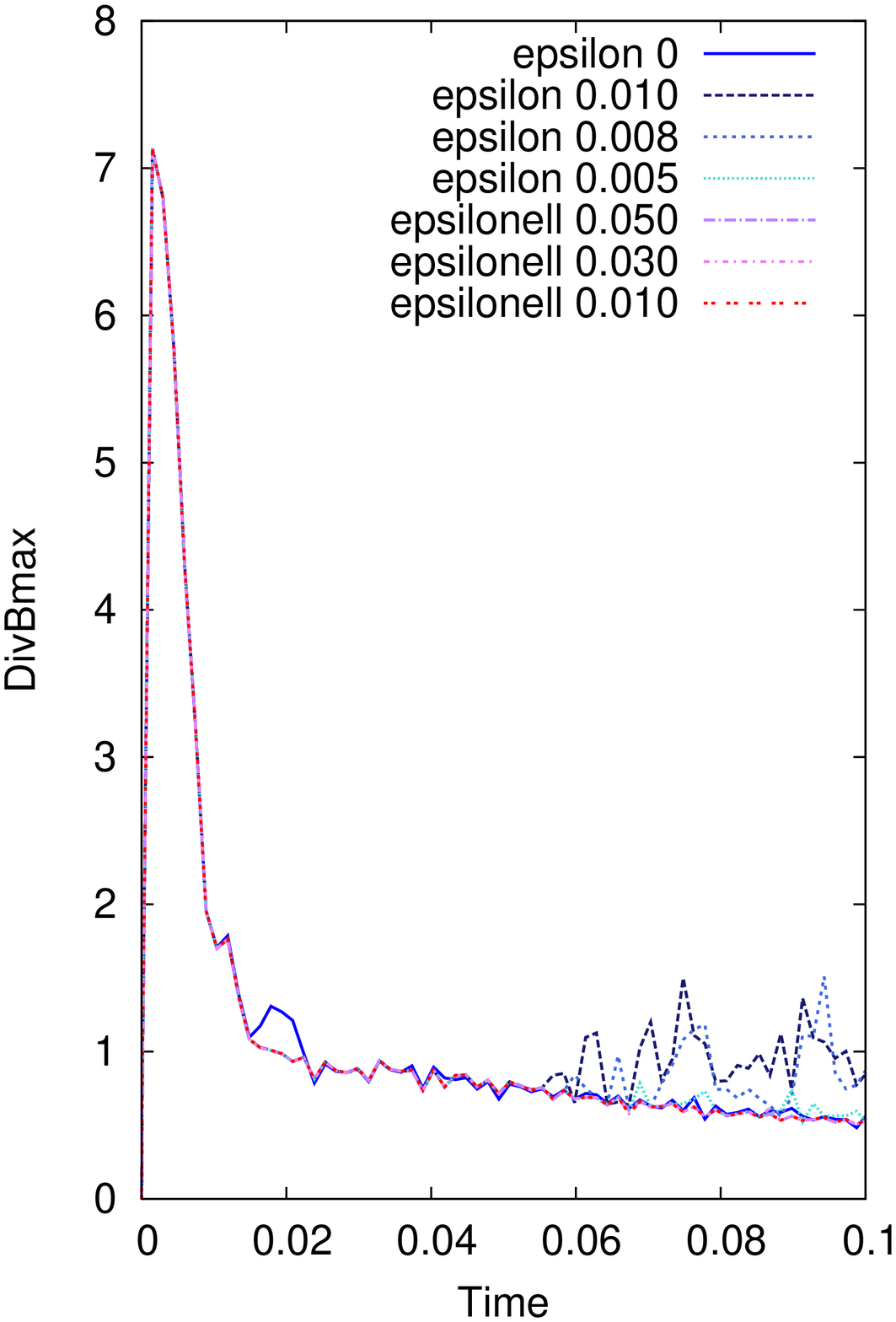}\\
 $(c)\,$MR, $\,L=10$ & $(d)\,$FV, $\,L=11$ \\
\includegraphics[width=0.35\linewidth]{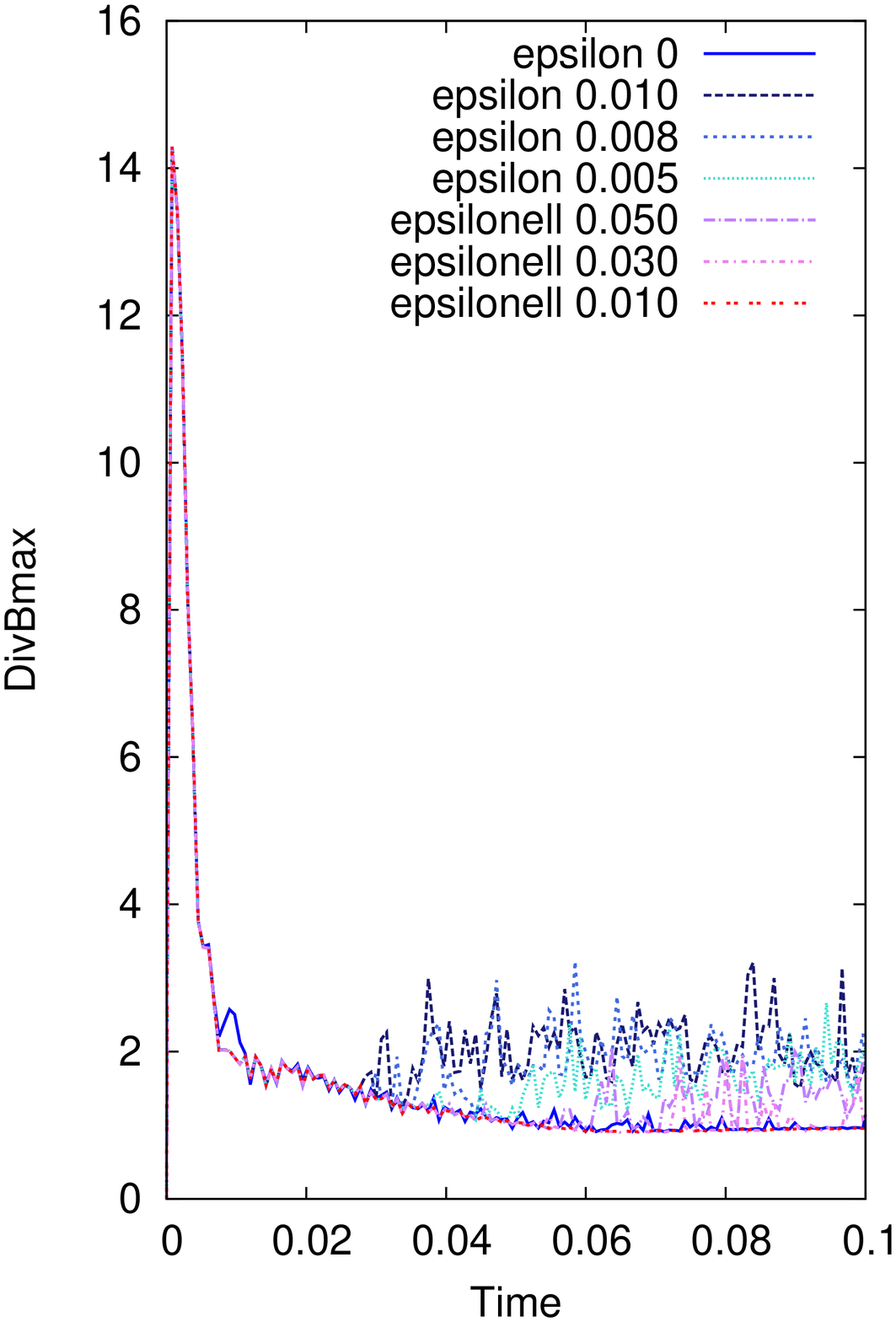} &
\includegraphics[width=0.35\linewidth]{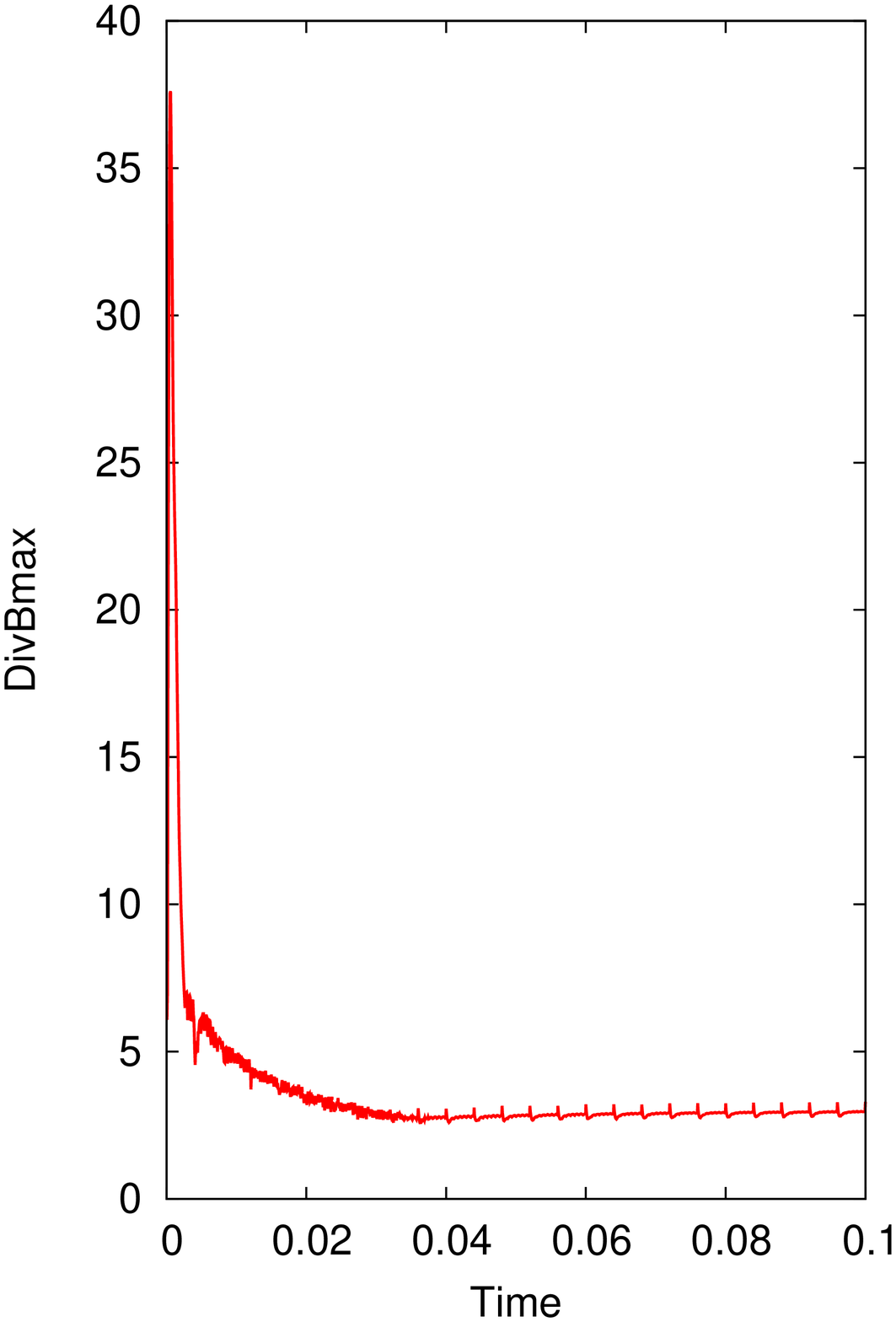}
\end{tabular}
\end{center}
\caption{The quantity $B_\mathrm{div}(t)$ over time for the  2D Riemann problem, with: $(a,b,c)$ GLM--MHD with mixed correction using the MR scheme with $\epsilon^\ell=\epsilon=0.010,\,0.008,\,0.005$ and $\epsilon^0=0,\,0.05,\,0.03,\,0.01$ for $L=8,\,9,\,10$; $(d)$ GLM--MHD with mixed correction using the FV scheme for the reference solution with $L=11$. 
}
\label{fig:2DRdivB}
\end{figure}

\begin{figure}[htb]
\psfrag{epsilon 0}{\tiny{$\epsilon=0.$}}
\psfrag{epsilon 0.005}{\tiny{$\epsilon=0.005$}}
\psfrag{epsilonj 0.05}{\tiny{$\epsilon^0=0.05$}}
\psfrag{DivBmax}{$B_\mathrm{div}(t)$}
\psfrag{dH/dt}{$\frac{\partial H}{\partial t}$}
\begin{center}
\begin{tabular}{ccc}
\includegraphics[width=0.35\linewidth]{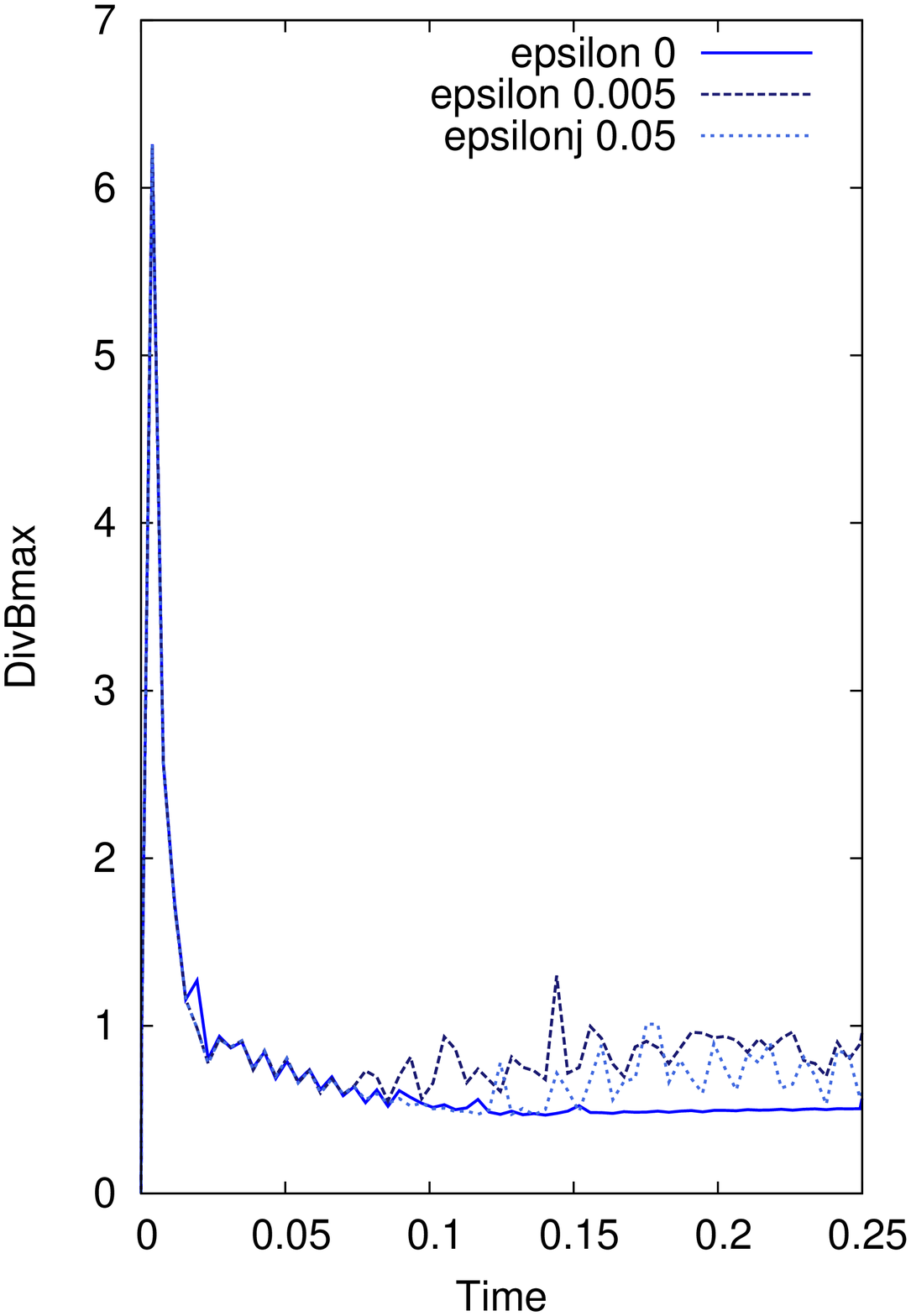}&
\includegraphics[width=0.35\linewidth]{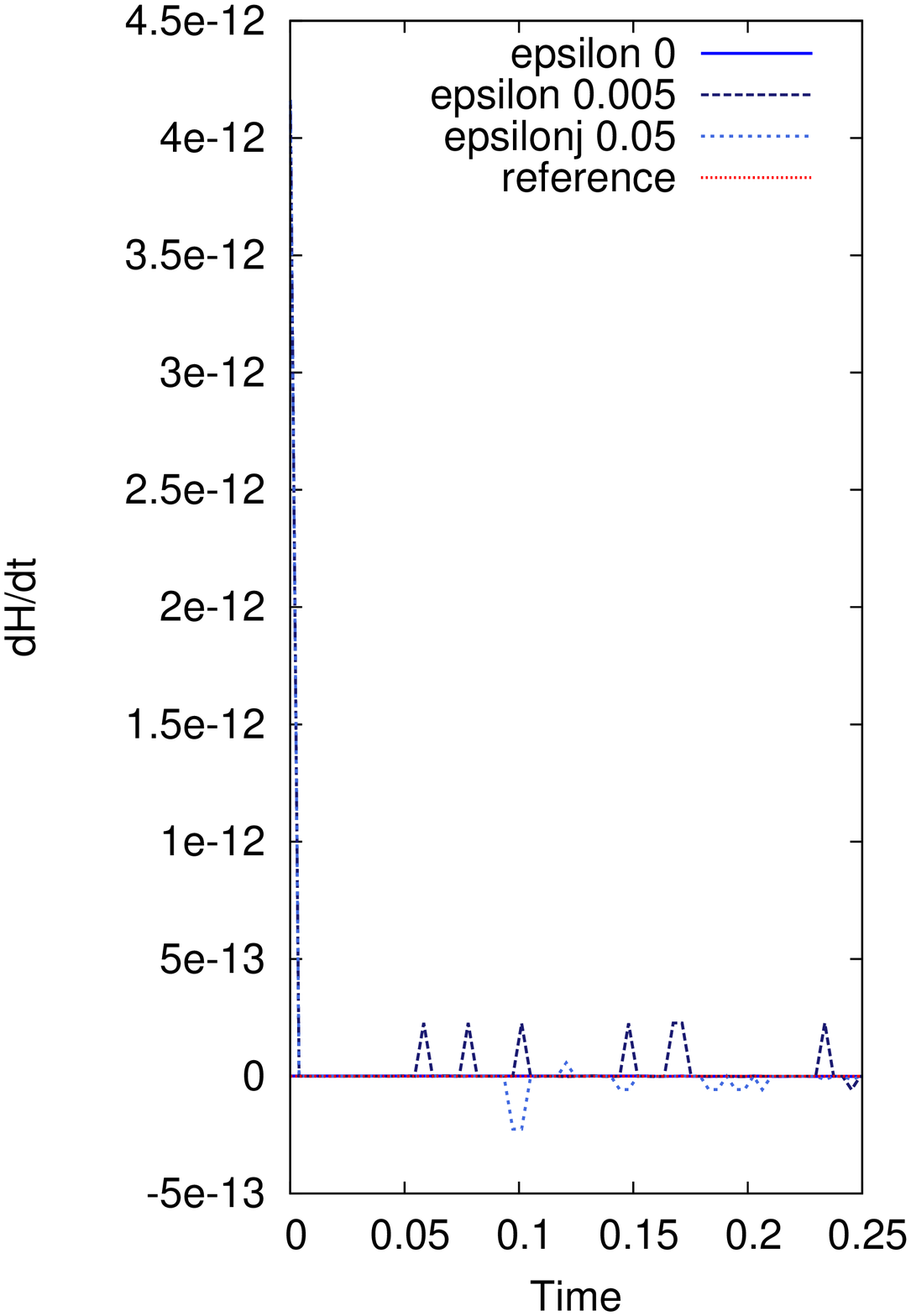}\\
\end{tabular}
\end{center} 
\caption{The quantities $B_\mathrm{div}(t)$ and time rate of change of magnetic helicity over time for the 2D Riemann problem, obtained with GLM--MHD with mixed correction MR scheme using $\epsilon^\ell=\epsilon=0,\,0.005$ and $\epsilon^0=0.05$ for $L=9$ and for reference solution.
}
\label{fig:2DRtescalar025}
\end{figure}

To check the time evolution of the divergence of the magnetic field, we consider the quantity
\[
B_\mathrm{div}(t):=\max\{|\nabla \cdot {\bf B}|:(x,y)\in[-1,1]^2\},
\]
where $\nabla \cdot \bf B$ is again evaluated using centered finite differences. Fig.~\ref{fig:2DRdivB} shows the time evolution of $B_\mathrm{div}(t)$ up to $t=0.1$ for the FV reference solution with $L=11$ (d) and three series of MR computations with $L=8,9,10$ (a, b, c) considering the following threshold, values $\epsilon = 0,\; 0.010,\; 0.008,\; 0.005$ and $\epsilon^0 = 0.050,\; 0.030,\; 0.010$. For the reference solution we observe a rapid decay of the initial value, around 37, during the first iterations, followed by a relaxation towards the value 3 which is reached at about 0.04. Afterwards, this value remains almost constant. For the MR computations we find that not only the initial but also the relaxation values of $B_\mathrm{div}(t)$ depend on the finest level $L$, and hence on the mesh size. For larger values of $L$ the divergence becomes larger but in all cases we find that after a certain time $B_\mathrm{div}(t)$ becomes constant or oscillates around a mean value. Using Harten's strategy with $\epsilon^0$ these oscillations almost disappear. In Fig.~\ref{fig:2DRtescalar025} we consider the evolution of $B_\mathrm{div}(t)$ for longer times, up to $t=0.25$, in MR cases with $L=9$ for $\epsilon=0$ and $0.005$, and $\epsilon^0=0.05$. After $t=0.1$ no oscillations can be observed for $\epsilon=0$, while for both $\epsilon^0 = 0.05$ and $\epsilon=0.005$ again some oscillations appear.

One main conclusion in analyzing $B_\mathrm{div}(t)$ for the different cases is that no growth in time can be observed, thus the divergence error seems to be controlled by the divergence cleaning, as discussed in \cite{Komissarov2007}.

Considering the conservative quantities \cite{yu2009note}, we compute the energy,
\[
\mathcal{E} = \int\int \left(|{\bf u}|^2 + |{\bf B}|^2\right)dxdy,
\]
and find the value $3.69$ at the initial time. At time $t=0.1$ we find for all FV solutions with $L=8,9$ and $10$ the value $3.48$. For the MR computations we obtain $3.46$, $3.47$ and $3.48$ for $L=8,9$ and $10$, respectively. These results are independent of the actual value of the threshold (ranging from $0.01$ down to $0$) and there is no significant influence if a fixed or level dependent value is used. This means that in all computations about $94\%$ of the energy is conserved. At a later time, $t=0.25$, we observe some decay, but still about $86\%$ of the energy is conserved.

The total magnetic helicity is also a conservative quantity of the ideal MHD equations \cite{bellan2006fundamentals} and we consider its time rate of change, defined as,
\[
\frac{\partial H}{\partial t} = a\int\int {\bf B}\cdot({\bf u}\times{\bf B})dxdy.
\]

As shown in Fig.~\ref{fig:2DRtescalar025}, right, the reference solution conserves perfectly the total magnetic helicity and $\partial H\ \partial t$ yields values close to the machine precision. For the three MR solutions there is an initial peak at about $4\cdot 10^{-12}$ which immediately decays to near zero machine precision, and remains zero for $\epsilon=0$. For the two others threshold values some intermittent spikes with amplitude below $2\cdot 10^{-13}$ are observed.


Table~\ref{tab:2DRCPUE} presents  a summary of the CPU time, memory compression, $D_c$ and $L_1^e(\rho)$ for all experiments at time $t=0.1$. 
For $\epsilon^{\ell}=\epsilon=0.005$ and $\epsilon^0=0.05$ the results are close, independent of the maximum level $L$.
However, the case  $\epsilon=0.005$ has slightly better CPU time and memory compression with respect to $L_1^e(\rho)$.
In these cases, for $L=10$, the CPU time are $7-14\%$ and the errors are approximately $10^{-2}$.
As expected, the error increases for a scale-independent threshold $\epsilon^{\ell}=\epsilon$ with $\epsilon$ being large, because it does not control well the error. 
However, as we decrease the value of $\epsilon$, the error becomes smaller. 
Thus, the choice of $\epsilon$ is an important ingredient. We can observe that if we choose a sufficiently small $\epsilon$, both strategies will have similar behavior.  However, we can optimize this process using Harten's strategy, which corresponds to a level dependent $\epsilon$.
\begin{table}[htb]
\caption{CPU time, memory, $D_c$, and density error $L_1^e(\rho)$ for the 2D Riemann problem computed with MR scheme using GLM--MHD with mixed correction and either with constant or level dependent threshold for $t=0.1$.}
  \label{tab:2DRCPUE}
\begin{center}
\begin{small}
\begin{tabular}{ccccp{0.1mm}cccc}
 \hline\\[-0.5mm]
  \multirow{2}{*}{$\boldsymbol{L=8}$} & \multicolumn{7}{c}{\textbf{MR}}& \multirow{1}{*}{\textbf{FV}} \\
  
   &  \multicolumn{3}{c}{$\epsilon^\ell=\epsilon$}&& \multicolumn{3}{c}{$\epsilon^0$} &\\
  	&   $0.01$ & $0.008$ & $0.005$ && $0.05$ & $0.03$ & $0.01$ &\\
\cline{2-4}\cline{6-8}
     CPU Time (\%)  & 22.74  & 23.47 & 24.55  && 26.71 & 27.80 & 30.33 & 100  \\  
     CPU Memory (\%) & 44.18 & 45.38 & 47.70 && 51.03 & 53.12 & 56.47 & 100 \\[0.5mm]
     $D_c$ (\%) & 29.74 & 30.67 & 32.50 && 34.94 & 36.60 & 39.28 & 100 \\[0.5mm]
     $L_1^e(\rho)\,\,\cdot 10^{-2}$ & 3.680 & 3.669 & 3.657 && 3.657 & 3.652 & 3.651 & 3.640\\[1mm]

  \hline\\[-0.5mm]
  \multirow{2}{*}{$\boldsymbol{L=9}$} & \multicolumn{7}{c}{\textbf{MR}}& \multirow{1}{*}{\textbf{FV}} \\

   &  \multicolumn{3}{c}{$\epsilon^\ell=\epsilon$}&& \multicolumn{3}{c}{$\epsilon^0$} &\\
  	&   $0.01$ & $0.008$ & $0.005$ && $0.05$ & $0.03$ & $0.01$ &\\
\cline{2-4}\cline{6-8}
	 CPU Time (\%)  & 13.63  & 14.66 & 15.91  && 17.67 & 19.00 & 20.46 & 100  \\    
     CPU Memory (\%) & 27.03 & 28.79 & 31.24 && 34.34 & 36.01 & 39.20 & 100 \\[0.5mm]
     $D_c$ (\%) & 17.70 & 18.97 & 21.01 && 23.51 & 24.92 & 27.42 & 100 \\[0.5mm]
     $L_1^e(\rho)\,\,\cdot 10^{-2}$ & 2086 & 2.039 & 1.981 & &1.974 & 1.958 & 1.953 & 1.9409\\[1mm]
    \hline\\[-0.5mm]

  \multirow{2}{*}{$\boldsymbol{L=10}$} & \multicolumn{7}{c}{\textbf{MR}}& \multirow{1}{*}{\textbf{FV}} \\
  
   &  \multicolumn{3}{c}{$\epsilon^\ell=\epsilon$}& &\multicolumn{3}{c}{$\epsilon^0$} &\\
  	&   $0.01$ & $0.008$ & $0.005$ & &$0.05$ & $0.03$ & $0.01$ &\\
\cline{2-4}\cline{6-8}
	 CPU Time (\%)  & 7.73  & 8.71 & 9.85 && 12.00 & 13.03 & 14.67 & 100  \\  
	 CPU Memory & 14.66 & 16.02 & 18.82 && 22.40 & 24.46 & 27.48 & 100 \\[0.5mm]
     $D_c$ (\%) & 9.25 & 10.07 & 12.01 && 14.66 & 1649 & 19.25 & 100 \\[0.5mm]	  
     $L_1^e(\rho)\,\,\cdot 10^{-2}$ & 1.090 & 1.031 & 0.932 && 0.905 & 0.895 & 0.851 & 0.841\\[1mm]  
    \hline
  \end{tabular}
\end{small}  
\end{center}
{\footnotesize NOTE: The results are computed with second order Runge-Kutta for the MR scheme. The CPU time for the GLM--MHD FV method is $277$ sec., $2326$ sec. and $314$ min., for $L=8,\,9$ and $10$,  at a Intel(R) Xeon(R) CPU E5620  \@ $2.40$GHz, CPU $1596$ MHz,  cache size $12288$ KB 
and $4$ cores. CPU time, memory and $D_c$ performances are computed with the corresponding non-adaptive FV solution using $L=8,\,9$ and $10$ scales on a uniform level. 
For the error, in all cases, we use a reference solution computed with a GLM--MHD FV scheme with $L=11$ for the same numerical scheme, implemented in the AMROC code \cite{Deiterdingetal:2009}.}
\end{table}

In Table~\ref{tab:2DRCPUE025} we show the CPU time, memory compression, $D_c$, and $L_1^e(\rho)$ for all experiments done at time $t=0.25$. 
We present the simulations for $\epsilon^{\ell}=\epsilon=0.005$ and $\epsilon^0=0.05$. 
The results at $t=0.25$ show that the MR approach does not introduce growing instabilities and it is possible to compute the solution for larger values of $t$.

\begin{table}[htb]
\caption{CPU time, memory, $D_c$, and density error $L_1^e(\rho)$ for the 2D Riemann problem simulated with the MR scheme using GLM--MHD with mixed correction and with constant or level dependent threshold for $t=0.25$}
  \label{tab:2DRCPUE025}
\begin{center}
\begin{small}
\begin{tabular}{ccp{0.1mm}cc}
%

 \hline\\[-0.5mm]
  \multirow{2}{*}{$\boldsymbol{L=9}$} & \multicolumn{3}{c}{\textbf{MR}}& \multirow{1}{*}{\textbf{FV}} \\
  
  	&   $\epsilon^\ell=0.005$ && $\epsilon^0=0.05$  &\\
\cline{2-4}
     CPU Time (\%)  & 18.79  && 22.61 &  100  \\  
     Memory (\%) & 38.12 && 45.25 & 100 \\[0.5mm]
     $D_c$ (\%) & 23.80 && 29.03 & 100 \\[0.5mm]
     $L_1^e(\rho)\,\,\cdot 10^{-2}$ & 3.887 && 3.826 & 3.694\\[1mm]
     
  \hline\\[-0.5mm]
%
  \end{tabular}
\end{small}  
\end{center}
\end{table}

\clearpage
\section{Conclusions and perspectives}
\label{sec:final}

Starting from the ideal MHD equations completed with generalized Lagrangian multipliers to
control the incompressibility of the magnetic field, we have developed an adaptive multiresolution method in two space dimensions on a Cartesian mesh with local refinement.
The space discretization is based on finite volumes with an HLLD numerical flux. For time integration an explicit Runge--Kutta scheme has been applied. 
To introduce a locally refined spatial mesh and also for local interpolation of the flux values Harten's cell average multiresolution analysis has been used.

To assess the efficiency and quality of this new adaptive scheme, we have considered a two-dimensional Riemann problem. We compared this numerical solution with adaptive MR results for different threshold values and two strategies of varying resolution levels. 
The numerical results show that the divergence cleaning can indeed work successfully with adaptive space discretizations.
The MR method with constant thresholding exhibits better CPU time performance but worse precision when compared to the level dependent threshold. 
The only drawback with respect to the level dependent threshold computations is that the number of cells 
on the adaptive mesh is increased.
We also observed that energy and time rate of change of magnetic helicity, both conserved quantities in the ideal MHD equations, remain indeed approximately conserved in our adaptive MR computations. 

In future work we plan to complete the adaptive method with time adaptivity using local and controlled
time stepping and to perform thus fully adaptive simulations in three space dimensions.
A second interesting direction is to move to non-ideal MHD, taking into account resistive effects
and finite values of the fluid viscosity to study the physics of reconnection of current sheets, especially in space physics applications.


\clearpage
\section*{Acknowledgements}

\begin{small}
M. O. D. and O. M. thankfully acknowledge financial support from MCTI/ FINEP /INFRINPE-1 (grant 01.12.0527.00), CAPES (grants $86/2010-29$), CNPq (grants $21224-6/2013-7, 483226/2011-4 , 306828/2010-3, 307511/2010-3,486165/2006-0, 305274/2009-0$), Ecole Centrale de Marseille (ECM), and 
FAPESP (grants $2012/06577-5, 2012/072812-2, 2007/07723-7$). 
A. G. thankfully acknowledges financial support for her Master, MCTI/INPE-PCI and PhD scholarship from CNPq (grants $132045/2010-9, 312479/2012-3, 141741/2013-9 $). 
K. S. thanks the ANR project SiCoMHD (ANR-Blanc 2011-045) for financial support.
%
%
We are grateful to Dominique Foug\`ere, Marie G. Dejean and Varlei E. Menconi (FAPESP grants 2008/09736-1 and MCTI/INPE-CNPq-PCI 312486/2012-0 and 455057/2013-5) for their helpful computational assistance.

\end{small}

\clearpage




\section*{Appendix HLLD Riemann Solver}\label{app:HLLD}
In the following solver, we consider the one-dimensional GLM-MHD equations in their primitive form
\begin{subequations}
\label{glm1d}
\begin{small}
\begin{eqnarray}
\displaystyle\frac{\partial\rho}{\partial t}
+\frac{\partial\rho u_x}{\partial x}&=&0,
\\
\displaystyle\frac{\partial E}{\partial t} +  \frac{\partial}{\partial x}\left[\left(E + p + \frac{{\bf B}^ 2}{2}\right){ u_x} - \left(u_x B_x + u_y B_y + u_z B_z \right) { B_x}\right]  &=& 0, 
\\
\displaystyle\frac{\partial \left(\rho { u_x}\right)}{\partial t} + \frac{\partial}{\partial x}\left[\rho u_x^2 + p +\left( \frac{{\bf B}^2}{2} \right)- B_x^2\right]  &=& { 0}, 
\\
\displaystyle\frac{\partial \left(\rho { u_y}\right)}{\partial t} + \frac{\partial}{\partial x}\left(\rho u_x u_y - B_xB_y\right)  &=& { 0}, 
\\
\displaystyle\frac{\partial \left(\rho { u_z}\right)}{\partial t} + \frac{\partial}{\partial x}\left(\rho u_x u_z - B_xB_z\right)  &=& { 0}, 
\\
\displaystyle\frac{\partial{ B_x}}{\partial t} + \frac{\partial \psi}{\partial x}&=&{ 0}, 
\\
\displaystyle\frac{\partial{ B_y}}{\partial t} + \frac{\partial}{\partial x}\left( u_xB_y - B_x u_y \right) &=&{ 0}, 
\\
\displaystyle\frac{\partial{ B_z}}{\partial t} + \frac{\partial}{\partial x}\left( u_xB_z - B_x u_z \right)   &=&{ 0}, 
\\
\frac{\partial\psi}{\partial t}+c_h^2\frac{\partial B_x}{\partial x}&=&-\frac{c_h^2}{c_p^2}\psi,
\end{eqnarray}
\end{small}
\end{subequations}
Considering the MHD system described above, we can obtain the Jacobian matrix. From the structure of this matrix one can verify that the equations of $B_x$ and $\psi$ can be decoupled from the remaining system and we can obtain the Jacobian matrix for the 1D MHD system \cite[p. 651-653]{Dedneretal:2002}. The eigenvalues of this matrix are $u_x$, $u_x\pm c_s$, $u_x\pm c_a$ and $u_x \pm c_f$, where $c_s,\,c_f$ are the slow and fast magneto-acoustic waves and $c_a$ is the Alfvén wave.

\begin{figure}[htb]
\psfrag{t}{$t$}
\psfrag{x}{$x$}
\psfrag{SM}{$S_M$}
\psfrag{SL}{$S_L$}
\psfrag{SLXT}{{\small $S_L\!=\!x/t$}}
\psfrag{SLS}{$S_L^\star$}
\psfrag{SRS}{$S_R^\star$}
\psfrag{QL}{$\textbf{Q}_L$}
\psfrag{QS}{$\textbf{Q}^\star$}
\psfrag{QLS}{$\textbf{Q}_L^\star$}
\psfrag{QLSS}{$\textbf{Q}_L^{\star\star}$}
\psfrag{QRSS}{$\textbf{Q}_R^{\star\star}$}
\psfrag{QRS}{$\textbf{Q}_R^\star$}
\psfrag{SR}{$S_R$}
\psfrag{SRXT}{{\small $S_R\!=\!x/t$}}
\psfrag{QR}{$\textbf{Q}_R$}
\begin{center}
\includegraphics[width=0.45\linewidth]{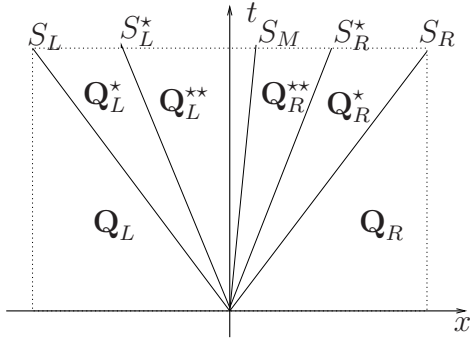}
\end{center}
\caption{Schematic of the Riemann fan structure with four intermediate states used in the HLLD flux.
Adapted from \cite{Kusano:2005}.}
\label{Fig:RiemannFanHLL}
\end{figure}

The Harten-Lax-van Leer-Discontinuities (HLLD) solver for MHD was firstly developed by Miyoshi and Kusano~\cite{Kusano:2005} and it can be considered as an extension of the Harten-Lax-van Leer (HLL) solver presented in~\cite{HLL:1983}. 
The HLLD solver is based on four intermediary states ${\bf Q}_L^{\star}$, ${\bf Q}_L^{\star\star}$, ${\bf Q}_R^{\star\star}$ and ${\bf Q}_R^{\star}$, divided by five waves $S_L$, $S_L^\star$, $S_M$, $S_R^\star$ and $S_R$, as illustrated in Fig.~\ref{Fig:RiemannFanHLL}. These waves are related to the entropy, fast and Alfvén waves. The HLLD numerical flux can resolve isolated discontinuities in the MHD system solution. This solver preserves positivity and it is more robust and efficient than the linearized solver, with an equally good resolution.

The states ${\bf Q}^\star$ and ${\bf Q}^{\star\star}$ for the GLM--MHD system are defined as \[{\bf Q}^\star_\alpha = (\rho^{\star}_\alpha, E^{\star}_\alpha, \rho^{\star}_\alpha {\bf u}^{\star}_\alpha,{\bf B}^{\star}_\alpha, \psi^{\star}_\alpha) \text{ and } {\bf Q}^{\star\star}_\alpha = (\rho^{\star\star}_\alpha, E^{\star\star}_\alpha, \rho^{\star\star}_\alpha {\bf u}^{\star\star}_\alpha,{\bf B}^{\star\star}_\alpha, \psi^{\star\star}_\alpha),\] with $\alpha$ denoting left ($L$) or right ($R$) states. 
In this approach,  we  compute the numerical flux for $\psi$ directly, then we consider $\psi^\star=\psi^{\star\star}=\psi$ here in the intermediary vector states, recalling that the HLLD is originally designed for MHD system, where the vector state $\mathbf{Q}$ has not the variable $\psi$ .

The numerical flux function is given by  
\begin{equation}
  	  	{\bf F}_{HLLD}=
  	  	\begin{cases}
  	  	\begin{array}{cl}
  	  	{\bf F}_L,		& \text{ if } S_L >0,\\
  	  	{\bf F}_L^\star,	& \text{ if } S_L \leq 0\leq S_L^\star,\\
  	  	{\bf F}_L^{\star\star},	& \text{ if } S_L^\star \leq 0\leq S_M,\\
  	  	{\bf F}_R^{\star\star},	& \text{ if } S_M \leq 0\leq S_R^\star,\\
  	  	{\bf F}_R^\star,	& \text{ if } S_R^\star \leq 0\leq S_R,\\
  	  	{\bf F}_R,		& \text{ if } S_R < 0.
  	  	\end{array}
  	  	\end{cases}.
  	  	\label{HLLD_flux}
  	  	\end{equation} 
The flux vectors ${\bf F}_L\,=\,{\bf F}({\bf Q}_L)$, ${\bf F}_R\,=\,{\bf F}({\bf Q}_R)$ are exact, while ${\bf F}_L^\star$, ${\bf F}_R^\star$ are approximate fluxes at intermediary states ${\bf Q}_L^{\star}$, ${\bf Q}_R^{\star}$, and ${\bf F}_L^{\star\star}$, ${\bf F}_R^{\star\star}$ are approximate fluxes at intermediary states ${\bf Q}_L^{\star\star}$, ${\bf Q}_R^{\star\star}$.

By the following process, we present the variables of the states ${\bf Q}_\alpha^{\star}$ and ${\bf Q}_\alpha^{\star\star}$, allowing us to compute the HLLD flux in the intermediary states
\begin{equation}
\label{HLLD2}
\begin{array}{l}
  	  	{\bf F}_\alpha^\star = {\bf F}_\alpha + S_\alpha\,({\bf Q}_\alpha^\star - {\bf Q}_\alpha), \\
  	  	{\bf F}_\alpha^{\star\star} = {\bf F}_\alpha + S_\alpha^\star\,{\bf Q}_\alpha^{\star\star} - (S_\alpha^\star - S_\alpha)\,{\bf Q}_\alpha^\star - S_\alpha\,{\bf Q}_\alpha,
\end{array}
\end{equation}  	  	
where $\alpha=R$ and $L$ denote right and left, respectively. 

The following description of the HLLD flux is related to the $x$ direction, considering
$B_x^\star=B_x^{\star\star}=B_x$.
In two-dimension, a similar expression can be obtained in the $y$ direction, considering $B_y^\star=B_y^{\star\star}=B_y$.

There are different possibilities to approximate the propagation speeds $S_\alpha$;  	  	
for instance, we use
\begin{equation}
S_L=\min(u_L,u_R)-\max(c_{f_L},c_{f_R}), \quad S_R=\max(u_L,u_R)+\max(c_{f_L},c_{f_R}),
\label{slopes}
\end{equation}
where $u_\alpha$ are the plasma velocities,  $c_{f_\alpha}$ are the magnetic acoustic waves \cite{Powell:1994}. 
The choice of $S_M$ is made to estimate the average normal velocity and it is given by 
	  \begin{equation}
	  S_M=\frac{(S_R-u_{x_R})\rho_R\; u_{x_R} - (S_L - u_{x_L})\rho_L\; u_{x_L} - p_{T_R} + p_{T_L}}{(S_R - u_{x_R})\rho_R - (S_L - u_{x_L})\rho_L}.
	  \label{eq:waveSM}
	  \end{equation}
The velocity is assumed to be constant over the Riemann fan, \textit{i.e.},
	  \begin{equation}
	  u_{x_L}^\star = u_{x_L}^{\star\star}=u_{x_R}^{\star\star}=u_{x_R}^\star=S_M.
	  \label{eq:velocityC}
	  \end{equation}
The total pressure  $p_{T} = p + \frac{{\bf B \cdot B}}{2}$ is kept constant, 
\begin{equation}
	  p_{T_L}^\star=p_{T_L}^{\star\star}=p_{T_R}^{\star\star}=p_{T_R}^\star=p_{T}^\star.
\label{eq:pressaoTotal}
\end{equation}
Under these conditions tangential and rotational discontinuities can be formed in the Riemann fan.

From the choice of $S_M$, the pressure $p_{T}^\star$ can be written as
 \begin{footnotesize}
\begin{equation}
\begin{split}
	  p_{T}^\star & =  \frac{(S_R - u_{x_R})\rho_R\; p_{T_L} - (S_L - u_{x_L})\rho_L\; p_{T_R}}{(S_R - u_{x_R})\rho_R - (S_L - u_{x_L})\rho_L} \\   & + \frac{\rho_L\; \rho_R (S_R - u_{x_R})(S_L - u_{x_L})(u_{x_R} - u_{x_L})}{(S_R - u_{x_R})\rho_R - (S_L - u_{x_L})\rho_L}.
	  \end{split}
	  \label{eq:pTotalstar}
	  \end{equation}
 \end{footnotesize}
 
Given $S_M$ and $p_T^\star$, the states ${\bf Q_\alpha^\star}=(\rho^\star_\alpha,p^\star_\alpha,u_{x_\alpha}^\star
	  ,u_{y_\alpha}^\star,u_{z_\alpha}^\star,B_{x_\alpha}^\star,
	  B_{y_\alpha}^\star,B_{z_\alpha}^\star)$ are bordered by the states $\bf Q_\alpha$ and they can be obtained from the jumps along $S_\alpha$, where $\alpha=L$ or $R$ represents the left or right state.
Therefore, one can derive the variables of the states ${\bf Q}_\alpha^\star$ as
	  \begin{subequations}
	  \begin{eqnarray}
	  \rho_\alpha^\star &=& \rho_\alpha\frac{S_\alpha- u_{x_\alpha}}{S_\alpha - S_M}, \\
	  u_{y_\alpha}^\star&=&u_{y_\alpha} - B_x B_{y_\alpha}\frac{S_M - u_{x_\alpha}}{\rho_\alpha (S_\alpha - u_{x_\alpha})(S_\alpha - S_M)  -B_x^2}, \\
	  u_{z_\alpha}^\star&=&u_{z_\alpha} - B_x B_{z_\alpha}\frac{S_M - u_{x_\alpha}}{\rho_\alpha (S_\alpha - u_{x_\alpha})(S_\alpha - S_M)  -B_x^2}, \\
	  B_{y_\alpha}^\star&=&B_{y_\alpha}\frac{\rho_\alpha(S_\alpha - u_{x_\alpha})^2 - B_x^2}{\rho_\alpha (S_\alpha - u_{x_\alpha})(S_\alpha - S_M)-B_x^2}, \\
	  B_{z_\alpha}^\star&=&B_{z_\alpha}\frac{\rho_\alpha(S_\alpha - u_{x_\alpha})^2 - B_x^2}{\rho_\alpha (S_\alpha - u_{x_\alpha})(S_\alpha - S_M)-B_x^2}.
	  \end{eqnarray}
	  \label{eq:rho_star}
	  \end{subequations}
Consequently, we can compute  $E_\alpha^\star$
	  \begin{equation}
	  E_\alpha^\star = \frac{(S_\alpha - u_{x_\alpha})E_\alpha - p_{T_\alpha}u_{x_\alpha} + p_T^\star S_M + B_x({\bf u}_\alpha\cdot{\bf B}_\alpha - {\bf u}_\alpha^\star \cdot{\bf B}_\alpha^\star)}{S_\alpha - S_M}.
	  \label{eq:energia_star}
	  \end{equation}

During the computations some operations as $0/0$ can appear when $S_M=u_{x_\alpha}$, $S_\alpha=u_{x_\alpha} \pm c_{f_\alpha}$, $B_{y_\alpha}=B_{z_\alpha}=0$ and $B_x^2\geq\gamma p_\alpha$. 
In these cases, we have to replace  $u_{y_\alpha}^\star=u_{y_\alpha}$, $u_{z_\alpha}^\star=u_{z_\alpha}$, and $B_{y_\alpha}^\star=B_{z_\alpha}^\star=0$.

Similarly, it is possible to obtain the equations related to the states  
\[{\bf Q_\alpha^{\star\star}}= (\rho^{\star\star}_\alpha,p^{\star\star}_\alpha,u_{x_\alpha}^{\star\star}
	  ,u_{y_\alpha}^{\star\star},u_{z_\alpha}^{\star\star},B_{x_\alpha}^{\star\star},
	  B_{y_\alpha}^{\star\star},B_{z_\alpha}^{\star\star}).\]
	  
Due to the relation described by Eq.~\ref{eq:velocityC}, starting with the jump condition of the continuity equation over an arbitrary value $S$, where $S_L<S<S_M$ or $S_M<S<S_R$, we have
	 \begin{equation}
	 \rho_\alpha^{\star\star} = \rho_\alpha^{\star}.
	 \label{eq:rho_starstar}
	 \end{equation}
The propagation velocities of the Alfv\'en waves in the intermediary states are estimated by
	 \begin{equation}
	 S_L^\star = S_M - \frac{|B_x|}{\sqrt{\rho_L^\star}}, \qquad S_R^\star = S_M - \frac{|B_x|}{\sqrt{\rho_R^\star}}.
	 \label{eq:alfven}
	 \end{equation}
Considering the jump conditions to the tangential components of the velocity and magnetic field over $S_M$, and  if $B_x\neq 0$, we can obtain the following relations
	 \begin{subequations}
	 \begin{eqnarray}
	 u_{y_L}^{\star\star}&=&u_{y_R}^{\star\star}\equiv u_{y}^{\star\star},\quad 
	 \label{eq:vystarstar}
	 u_{z_L}^{\star\star}=u_{z_R}^{\star\star}\equiv u_{z}^{\star\star},\\
	 \label{eq:vzstarstar}
	 B_{y_L}^{\star\star}&=&B_{y_R}^{\star\star}\equiv B_{y}^{\star\star},\quad
	 B_{z_L}^{\star\star} =B_{z_R}^{\star\star}\equiv B_{z}^{\star\star}.
	 \end{eqnarray}
	 \label{eq:Ustarstar}
	 \end{subequations}
If $B_x=0$, it is impossible to calculate the remaining variables of the states ${\bf Q_\alpha^{\star\star}}$.
Replacing  Eqs.~\ref{eq:rho_starstar} -- \ref{eq:Ustarstar} into the integral conservation laws over the Riemann fan, we can derive the variables
\begin{subequations}
\begin{eqnarray}
u_{y}^{\star\star} &=& \frac{ u_{y_L}^{\star} + \sqrt{\rho_R^\star} u_{y_R}^{\star} + (B_{y_R}^\star - B_{y_L}^\star)\text{sign}(B_x)}{\sqrt{\rho_L^\star} + \sqrt{\rho_R^\star}}, \\
u_{z}^{\star\star} &=& \frac{\sqrt{\rho_L^\star} u_{z_L}^{\star} + \sqrt{\rho_R^\star} u_{z_R}^{\star} + (B_{z_R}^\star - B_{z_L}^\star)\text{sign}(B_x)}{\sqrt{\rho_L^\star} + \sqrt{\rho_R^\star}}, \\
B_{y}^{\star\star} &=& \frac{\sqrt{\rho_L^\star} B_{y_R}^{\star} + \sqrt{\rho_R^\star} B_{y_L}^{\star} + \sqrt{\rho_L^\star \rho_R^\star}(u_{y_R}^\star - u_{y_L}^\star)\text{sign}(B_x)}{\sqrt{\rho_L^\star} + \sqrt{\rho_R^\star}}, \\
B_{z}^{\star\star} &=& \frac{\sqrt{\rho_L^\star} B_{z_R}^{\star} + \sqrt{\rho_R^\star} B_{z_L}^{\star} + \sqrt{\rho_L^\star \rho_R^\star}(u_{z_R}^\star - u_{z_L}^\star)\text{sign}(B_x)}{\sqrt{\rho_L^\star} + \sqrt{\rho_R^\star}},
\end{eqnarray}
\end{subequations}
where $\text{sign}(B_x)$ is $1$ for $B_x>0$, and $-1$ for $B_x<0$. 
Consequently, the equation of the energy  in $\bf Q^{\star\star}$ is given by
	\begin{equation}
	E_\alpha^{\star\star} = E_\alpha^{\star} \, \mp \,  \sqrt{\rho_\alpha^\star} \left({\bf u}_\alpha^\star\cdot {\bf B}_\alpha^\star - {\bf u}_\alpha^{\star\star}\cdot {\bf B}_\alpha^{\star\star}\right) \; \text{sign}(B_x).
	\label{eq:Estarstar}
	\end{equation}
The same procedure is done for the $y$ direction.
\clearpage
\bibliographystyle{plain}


\end{document}